
%
%
%
\def\unredoffs{} \def\redoffs{\voffset=-.31truein\hoffset=-.59truein}
\def\speclscape{\special{ps: landscape}}
%
%
%
%
\newbox\leftpage \newdimen\fullhsize \newdimen\hstitle \newdimen\hsbody
\tolerance=1000\hfuzz=2pt
\catcode`\@=11 
%
\ifx\answ\bigans\message{(This will come out unreduced.}
\magnification=1200\unredoffs\baselineskip=16pt plus 2pt minus 1pt
\hsbody=\hsize \hstitle=\hsize 
\else\message{(This will be reduced.} \let\l@r=L
\magnification=1000\baselineskip=16pt plus 2pt minus 1pt \vsize=7truein
\redoffs \hstitle=8truein\hsbody=4.75truein\fullhsize=10truein\hsize=\hsbody
\output={\ifnum\pageno=0 
  \shipout\vbox{\speclscape{\hsize\fullhsize\makeheadline}
    \hbox to \fullhsize{\hfill\pagebody\hfill}}\advancepageno
  \else
  \almostshipout{\leftline{\vbox{\pagebody\makefootline}}}\advancepageno 
  \fi}
\def\almostshipout#1{\if L\l@r \count1=1 \message{[\the\count0.\the\count1]}
      \global\setbox\leftpage=#1 \global\let\l@r=R
 \else \count1=2
  \shipout\vbox{\speclscape{\hsize\fullhsize\makeheadline}
      \hbox to\fullhsize{\box\leftpage\hfil#1}}  \global\let\l@r=L\fi}
\fi
%
\newcount\yearltd\yearltd=\year\advance\yearltd by -1900

\def\Title#1#2{\nopagenumbers\abstractfont\hsize=\hstitle\rightline{#1}%
\vskip 1in\centerline{\titlefont #2}\abstractfont\vskip .5in\pageno=0}
\def\Date#1{\vfill\leftline{#1}\tenpoint\supereject\global\hsize=\hsbody%
\footline={\hss\tenrm\folio\hss}}
%

\def\draftmode{\message{ DRAFTMODE }\def\draftdate{{\rm preliminary draft:
\number\month/\number\day/\number\yearltd\ \ \hourmin}}%
\headline={\hfil\draftdate}\writelabels\baselineskip=20pt plus 2pt minus 2pt
 {\count255=\time\divide\count255 by 60 \xdef\hourmin{\number\count255}
  \multiply\count255 by-60\advance\count255 by\time
  \xdef\hourmin{\hourmin:\ifnum\count255<10 0\fi\the\count255}}}
\def\nolabels{\def\wrlabeL##1{}\def\eqlabeL##1{}\def\reflabeL##1{}}
\def\writelabels{\def\wrlabeL##1{\leavevmode\vadjust{\rlap{\smash%
{\line{{\escapechar=` \hfill\rlap{\sevenrm\hskip.03in\string##1}}}}}}}%
\def\eqlabeL##1{{\escapechar-1\rlap{\sevenrm\hskip.05in\string##1}}}%
\def\reflabeL##1{\noexpand\llap{\noexpand\sevenrm\string\string\string##1}}}
\nolabels
%
\global\newcount\secno \global\secno=0
\global\newcount\meqno \global\meqno=1
\def\newsec#1{\global\advance\secno by1\message{(\the\secno. #1)}
\global\subsecno=0\eqnres@t\noindent{\bf\the\secno. #1}
\writetoca{{\secsym} {#1}}\par\nobreak\medskip\nobreak}
\def\eqnres@t{\xdef\secsym{\the\secno.}\global\meqno=1\bigbreak\bigskip}
\def\sequentialequations{\def\eqnres@t{\bigbreak}}\xdef\secsym{}
\global\newcount\subsecno \global\subsecno=0
\def\subsec#1{\global\advance\subsecno by1\message{(\secsym\the\subsecno. #1)}
\ifnum\lastpenalty>9000\else\bigbreak\fi
\noindent{\it\secsym\the\subsecno. #1}\writetoca{\string\quad 
{\secsym\the\subsecno.} {#1}}\par\nobreak\medskip\nobreak}
\def\appendix#1#2{\global\meqno=1\global\subsecno=0\xdef\secsym{\hbox{#1.}}
\bigbreak\bigskip\noindent{\bf Appendix #1. #2}\message{(#1. #2)}
\writetoca{Appendix {#1.} {#2}}\par\nobreak\medskip\nobreak}
%
%
\def\eqnn#1{\xdef #1{(\secsym\the\meqno)}\writedef{#1\leftbracket#1}%
\global\advance\meqno by1\wrlabeL#1}
\def\eqna#1{\xdef #1##1{\hbox{$(\secsym\the\meqno##1)$}}
\writedef{#1\numbersign1\leftbracket#1{\numbersign1}}%
\global\advance\meqno by1\wrlabeL{#1$\{\}$}}
\def\eqn#1#2{\xdef #1{(\secsym\the\meqno)}\writedef{#1\leftbracket#1}%
\global\advance\meqno by1$$#2\eqno#1\eqlabeL#1$$}
%
\newskip\footskip\footskip14pt plus 1pt minus 1pt 
\def\footnotefont{\ninepoint}\def\f@t#1{\footnotefont #1\@foot}
\def\f@@t{\baselineskip\footskip\bgroup\footnotefont\aftergroup\@foot\let\next}
\setbox\strutbox=\hbox{\vrule height9.5pt depth4.5pt width0pt}
\global\newcount\ftno \global\ftno=0
\def\foot{\global\advance\ftno by1\footnote{$^{\the\ftno}$}}
%
\newwrite\ftfile   
\def\footend{\def\foot{\global\advance\ftno by1\chardef\wfile=\ftfile
$^{\the\ftno}$\ifnum\ftno=1\immediate\openout\ftfile=foots.tmp\fi%
\immediate\write\ftfile{\noexpand\smallskip%
\noexpand\item{f\the\ftno:\ }\pctsign}\findarg}%
\def\footatend{\vfill\eject\immediate\closeout\ftfile{\parindent=20pt
\centerline{\bf Footnotes}\nobreak\bigskip\input foots.tmp }}}
\def\footatend{}
%
%
\global\newcount\refno \global\refno=1
\newwrite\rfile
%
\def\ref{\nref}
\def\nref#1{\xdef#1{[\the\refno]}\writedef{#1\leftbracket#1}%
\ifnum\refno=1\immediate\openout\rfile=refs.tmp\fi
\global\advance\refno by1\chardef\wfile=\rfile\immediate
\write\rfile{\noexpand\item{#1\ }\reflabeL{#1\hskip.31in}\pctsign}\findarg}
\def\findarg#1#{\begingroup\obeylines\newlinechar=`\^^M\pass@rg}
{\obeylines\gdef\pass@rg#1{\writ@line\relax #1^^M\hbox{}^^M}%
\gdef\writ@line#1^^M{\expandafter\toks0\expandafter{\striprel@x #1}%
\edef\next{\the\toks0}\ifx\next\em@rk\let\next=\endgroup\else\ifx\next\empty%
\else\immediate\write\wfile{\the\toks0}\fi\let\next=\writ@line\fi\next\relax}}
\def\striprel@x#1{} \def\em@rk{\hbox{}} 
\def\lref{\begingroup\obeylines\lr@f}
\def\lr@f#1#2{\gdef#1{\ref#1{#2}}\endgroup\unskip}

\def\addref#1{\immediate\write\rfile{\noexpand\item{}#1}} 
\def\startrefs#1{\immediate\openout\rfile=refs.tmp\refno=#1}
\def\xref{\expandafter\xr@f}\def\xr@f[#1]{#1}
\def\refs#1{\count255=1[\r@fs #1{\hbox{}}]}
\def\r@fs#1{\ifx\und@fined#1\message{reflabel \string#1 is undefined.}%
\nref#1{need to supply reference \string#1.}\fi%
\vphantom{\hphantom{#1}}\edef\next{#1}\ifx\next\em@rk\def\next{}%
\else\ifx\next#1\ifodd\count255\relax\xref#1\count255=0\fi%
\else#1\count255=1\fi\let\next=\r@fs\fi\next}
%

%
\newwrite\ffile\global\newcount\figno \global\figno=1
\def\fig{fig.~\the\figno\nfig}
\def\nfig#1{\xdef#1{fig.~\the\figno}%
\writedef{#1\leftbracket fig.\noexpand~\the\figno}%
\ifnum\figno=1\immediate\openout\ffile=figs.tmp\fi\chardef\wfile=\ffile%
\immediate\write\ffile{\noexpand\medskip\noexpand\item{Fig.\ \the\figno. }
\reflabeL{#1\hskip.55in}\pctsign}\global\advance\figno by1\findarg}
\def\xfig{\expandafter\xf@g}\def\xf@g fig.\penalty\@M\ {}
\def\figs#1{figs.~\f@gs #1{\hbox{}}}
\def\f@gs#1{\edef\next{#1}\ifx\next\em@rk\def\next{}\else
\ifx\next#1\xfig #1\else#1\fi\let\next=\f@gs\fi\next}
\newwrite\lfile
{\escapechar-1\xdef\pctsign{\string\%}\xdef\leftbracket{\string\{}
\xdef\rightbracket{\string\}}\xdef\numbersign{\string\#}}

\def\writestop{\def\writestoppt{\immediate\write\lfile{\string\pageno%
\the\pageno\string\startrefs\leftbracket\the\refno\rightbracket%
\string\def\string\secsym\leftbracket\secsym\rightbracket%
\string\secno\the\secno\string\meqno\the\meqno}\immediate\closeout\lfile}}
\def\writestoppt{}\def\writedef#1{}
\def\seclab#1{\xdef #1{\the\secno}\writedef{#1\leftbracket#1}\wrlabeL{#1=#1}}
\def\subseclab#1{\xdef #1{\secsym\the\subsecno}%
\writedef{#1\leftbracket#1}\wrlabeL{#1=#1}}
\newwrite\tfile \def\writetoca#1{}
\def\leaderfill{\leaders\hbox to 1em{\hss.\hss}\hfill}
\def\writetoc{\immediate\openout\tfile=toc.tmp 
   \def\writetoca##1{{\edef\next{\write\tfile{\noindent ##1 
   \string\leaderfill {\noexpand\number\pageno} \par}}\next}}}
%
%
%

\catcode`\@=12 
%
\edef\tfontsize{\ifx\answ\bigans scaled\magstep3\else scaled\magstep4\fi}
\font\titlerm=cmr10 \tfontsize \font\titlerms=cmr7 \tfontsize
\font\titlermss=cmr5 \tfontsize \font\titlei=cmmi10 \tfontsize
\font\titleis=cmmi7 \tfontsize \font\titleiss=cmmi5 \tfontsize
\font\titlesy=cmsy10 \tfontsize \font\titlesys=cmsy7 \tfontsize
\font\titlesyss=cmsy5 \tfontsize \font\titleit=cmti10 \tfontsize
\skewchar\titlei='177 \skewchar\titleis='177 \skewchar\titleiss='177
\skewchar\titlesy='60 \skewchar\titlesys='60 \skewchar\titlesyss='60
\def\titlefont{\def\rm{\fam0\titlerm}
\textfont0=\titlerm \scriptfont0=\titlerms \scriptscriptfont0=\titlermss
\textfont1=\titlei \scriptfont1=\titleis \scriptscriptfont1=\titleiss
\textfont2=\titlesy \scriptfont2=\titlesys \scriptscriptfont2=\titlesyss
\textfont\itfam=\titleit \def\it{\fam\itfam\titleit}\rm}
 \ifx\answ\bigans\else scaled\magstep1\fi
\ifx\answ\bigans\def\abstractfont{\tenpoint}\else
\font\abssl=cmsl10 scaled \magstep1
\font\absrm=cmr10 scaled\magstep1 \font\absrms=cmr7 scaled\magstep1
\font\absrmss=cmr5 scaled\magstep1 \font\absi=cmmi10 scaled\magstep1
\font\absis=cmmi7 scaled\magstep1 \font\absiss=cmmi5 scaled\magstep1
\font\abssy=cmsy10 scaled\magstep1 \font\abssys=cmsy7 scaled\magstep1
\font\abssyss=cmsy5 scaled\magstep1 \font\absbf=cmbx10 scaled\magstep1
\skewchar\absi='177 \skewchar\absis='177 \skewchar\absiss='177
\skewchar\abssy='60 \skewchar\abssys='60 \skewchar\abssyss='60
\def\abstractfont{\def\rm{\fam0\absrm}
\textfont0=\absrm \scriptfont0=\absrms \scriptscriptfont0=\absrmss
\textfont1=\absi \scriptfont1=\absis \scriptscriptfont1=\absiss
\textfont2=\abssy \scriptfont2=\abssys \scriptscriptfont2=\abssyss
\textfont\itfam=\bigit \def\it{\fam\itfam\bigit}\def\footnotefont{\tenpoint}%
\textfont\slfam=\abssl \def\sl{\fam\slfam\abssl}%
\textfont\bffam=\absbf \def\bf{\fam\bffam\absbf}\rm}\fi
\def\tenpoint{\def\rm{\fam0\tenrm}
\textfont0=\tenrm \scriptfont0=\sevenrm \scriptscriptfont0=\fiverm
\textfont1=\teni  \scriptfont1=\seveni  \scriptscriptfont1=\fivei
\textfont2=\tensy \scriptfont2=\sevensy \scriptscriptfont2=\fivesy
\textfont\itfam=\tenit \def\it{\fam\itfam\tenit}\def\footnotefont{\ninepoint}%
\textfont\bffam=\tenbf \def\bf{\fam\bffam\tenbf}\def\sl{\fam\slfam\tensl}\rm}
\font\ninerm=cmr9 \font\sixrm=cmr6 \font\ninei=cmmi9 \font\sixi=cmmi6 
\font\ninesy=cmsy9 \font\sixsy=cmsy6 \font\ninebf=cmbx9 
\font\nineit=cmti9 \font\ninesl=cmsl9 \skewchar\ninei='177
\skewchar\sixi='177 \skewchar\ninesy='60 \skewchar\sixsy='60 
\def\ninepoint{\def\rm{\fam0\ninerm}
\textfont0=\ninerm \scriptfont0=\sixrm \scriptscriptfont0=\fiverm
\textfont1=\ninei \scriptfont1=\sixi \scriptscriptfont1=\fivei
\textfont2=\ninesy \scriptfont2=\sixsy \scriptscriptfont2=\fivesy
\textfont\itfam=\ninei \def\it{\fam\itfam\nineit}\def\sl{\fam\slfam\ninesl}%
\textfont\bffam=\ninebf \def\bf{\fam\bffam\ninebf}\rm} 
%
%

\hyphenation{anom-aly anom-alies coun-ter-term coun-ter-terms}
\def\inv{^{\raise.15ex\hbox{${\scriptscriptstyle -}$}\kern-.05em 1}}

\def\Dsl{\,\raise.15ex\hbox{/}\mkern-13.5mu D} 
\def\dsl{\raise.15ex\hbox{/}\kern-.57em\partial}

\font\bigit=cmti10 scaled \magstep1
\def\lspace{\ifx\answ\bigans{}\else\qquad\fi}
\def\lbspace{\ifx\answ\bigans{}\else\hskip-.2in\fi} 
\def\boxeqn#1{\vcenter{\vbox{\hrule\hbox{\vrule\kern3pt\vbox{\kern3pt
	\hbox{${\displaystyle #1}$}\kern3pt}\kern3pt\vrule}\hrule}}}
\def\mbox#1#2{\vcenter{\hrule \hbox{\vrule height#2in
		\kern#1in \vrule} \hrule}}  
%

\def\darr#1{\raise1.5ex\hbox{$\leftrightarrow$}\mkern-16.5mu #1}

\def\roughly#1{\raise.3ex\hbox{$#1$\kern-.75em\lower1ex\hbox{$\sim$}}}

\def\frac#1#2{{#1\over#2}}

\def\journal#1&#2(#3){\unskip, #1~\bf #2 \rm(19#3) }
\def\andjournal#1&#2(#3){\sl #1~\bf #2 \rm (19#3) }

\def\bra#1{\left\langle #1\right|}
\def\ket#1{\left| #1\right\rangle}
\def\det{{\rm det}}

\catcode`\@=11\def\slash#1{\mathord{\mathpalette\c@ncel{#1}}}
\overfullrule=0pt
\def\steepslash{\c@ncel}
\def\frac#1#2{{#1\over #2}}

\def\:{\!:\!}
\def\inbar{\,\vrule height1.5ex width.4pt depth0pt}
\def\IQ{\relax\,\hbox{$\inbar\kern-.3em{\rm Q}$}}
\def\IB{\relax{\rm I\kern-.18em B}}
\def\IC{\relax\hbox{$\inbar\kern-.3em{\rm C}$}}
\def\IP{\relax{\rm I\kern-.18em P}}
\def\IR{\relax{\rm I\kern-.18em R}}
\def\ZZ{\relax\ifmmode\mathchoice
{\hbox{Z\kern-.4em Z}}{\hbox{Z\kern-.4em Z}}
{\lower.9pt\hbox{Z\kern-.4em Z}}
{\lower1.2pt\hbox{Z\kern-.4em Z}}\else{Z\kern-.4em Z}\fi}

\catcode`\@=12

\def\npb#1(#2)#3{{ Nucl. Phys. }{B#1} (#2) #3}
\def\plb#1(#2)#3{{ Phys. Lett. }{#1B} (#2) #3}
\def\pla#1(#2)#3{{ Phys. Lett. }{#1A} (#2) #3}
\def\prl#1(#2)#3{{ Phys. Rev. Lett. }{#1} (#2) #3}
\def\mpla#1(#2)#3{{ Mod. Phys. Lett. }{A#1} (#2) #3}
\def\ijmpa#1(#2)#3{{ Int. J. Mod. Phys. }{A#1} (#2) #3}
\def\cmp#1(#2)#3{{ Comm. Math. Phys. }{#1} (#2) #3}
\def\cqg#1(#2)#3{{ Class. Quantum Grav. }{#1} (#2) #3}
\def\jmp#1(#2)#3{{ J. Math. Phys. }{#1} (#2) #3}
\def\anp#1(#2)#3{{ Ann. Phys. }{#1} (#2) #3}
\def\prd#1(#2)#3{{ Phys. Rev. } {D{#1}} (#2) #3}
\def\ptp#1(#2)#3{{ Progr. Theor. Phys. }{#1} (#2) #3}
\def\aom#1(#2)#3{{ Ann. Math. }{#1} (#2) #3}

\def\br{\buildrel}
\def\bra{\langle}
\def\ket{\rangle}

\def\C{{\bf C}}

\def\P{{\bf P}}
\def\Q{{\bf Q}}

\def\Z{{\bf Z}}
\def\cA{{\cal A}}

\def\cF{{\cal F}}

\def\cL{{\cal L}}
\def\cM{{\cal M}}

\def\cO{{\cal O}}

\def\cR{{\cal R}}
\def\cS{{\cal S}}
\def\cT{{\cal T}}
\def\cU{{\cal U}}

\def\cX{{\cal X}}

\def\cicy#1(#2|#3)#4{\left(\matrix{#2}\right|\!\!
                     \left|\matrix{#3}\right)^{{#4}}_{#1}}

\def\ra{\rightarrow}

\def\bs{\bigskip}

\def\Box{{\,\lower0.9pt\vbox{\hrule 
\hbox{\vrule height 0.2 cm \hskip 0.2 cm  
\vrule height 0.2 cm}\hrule}\,}}

\global\newcount\thmno \global\thmno=0
\def\definition#1{\global\advance\thmno by1
\bigskip\noindent{\bf Definition \secsym\the\thmno. }{\it #1}
\par\nobreak\medskip\nobreak}
\def\question#1{\global\advance\thmno by1
\bigskip\noindent{\bf Question \secsym\the\thmno. }{\it #1}
\par\nobreak\medskip\nobreak}
\def\theorem#1{\global\advance\thmno by1
\bigskip\noindent{\bf Theorem \secsym\the\thmno. }{\it #1}
\par\nobreak\medskip\nobreak}
\def\proposition#1{\global\advance\thmno by1
\bigskip\noindent{\bf Proposition \secsym\the\thmno. }{\it #1}
\par\nobreak\medskip\nobreak}
\def\corollary#1{\global\advance\thmno by1
\bigskip\noindent{\bf Corollary \secsym\the\thmno. }{\it #1}
\par\nobreak\medskip\nobreak}
\def\lemma#1{\global\advance\thmno by1
\bigskip\noindent{\bf Lemma \secsym\the\thmno. }{\it #1}
\par\nobreak\medskip\nobreak}
\def\conjecture#1{\global\advance\thmno by1
\bigskip\noindent{\bf Conjecture \secsym\the\thmno. }{\it #1}
\par\nobreak\medskip\nobreak}
\def\exercise#1{\global\advance\thmno by1
\bigskip\noindent{\bf Exercise \secsym\the\thmno. }{\it #1}
\par\nobreak\medskip\nobreak}
\def\remark#1{\global\advance\thmno by1
\bigskip\noindent{\bf Remark \secsym\the\thmno. }{\it #1}
\par\nobreak\medskip\nobreak}
\def\problem#1{\global\advance\thmno by1
\bigskip\noindent{\bf Problem \secsym\the\thmno. }{\it #1}
\par\nobreak\medskip\nobreak}
\def\others#1#2{\global\advance\thmno by1
\bigskip\noindent{\bf #1 \secsym\the\thmno. }{\it #2}
\par\nobreak\medskip\nobreak}
\def\proof{\noindent Proof: }

\def\thmlab#1{\xdef #1{\secsym\the\thmno}\writedef{#1\leftbracket#1}\wrlabeL{#1=#1}}
%
%
\def\newsec#1{\global\advance\secno by1\message{(\the\secno. #1)}
\global\subsecno=0\thmno=0\eqnres@t\noindent{\bf\the\secno. #1}
\writetoca{{\secsym} {#1}}\par\nobreak\medskip\nobreak}
\def\eqnres@t{\xdef\secsym{\the\secno.}\global\meqno=1\bigbreak\bigskip}
\def\sequentialequations{\def\eqnres@t{\bigbreak}}\xdef\secsym{}
%

%

\ref\AP{C. Allday, V. Puppe,{\it Cohomology methods in transformation groups},
Cambridge University Press, 1993.}
\ref\Atiyah{M. Atiyah, {\it Convexity and commuting Hamiltonians},
Bull. London Math. Soc. 14 (1982) 1-15.}
\ref\AB{M. Atiyah and R. Bott, {\it The moment map and
equivariant cohomology}, Topology 23 (1984) 1-28.} 
\ref\Batyrev{V. Batyrev, {\it Dual polyhedra and mirror symmetry
for Calabi-Yau hypersurfaces in toric varieties},
J. Alg. Geom. 3 (1994) 493-535.}
\ref\BB{V. Batyrev and L. Borisov, {\it 
On Calabi-Yau complete intersections in toric varieties},
alg-geom/9412017.}
\ref\BatyrevStraten{V. Batyrev and D. van Straten,
{\it Generalized hypergeometric functions and rational
curves on Calabi-Yau complete intersections in
toric varieties}, Comm. Math. Phys. 168 (1995) 495-533.}
\ref\BehrendFentachi{K. Behrend and B. Fantechi, 
{\it The intrinsic normal cone}, Invent. Math. 128 (1997) 45-88.}
\ref\BKK{P. Berglund, S. Katz and A. Klemm,
{\it Mirror Symmetry and the Moduli Space for Generic
Hypersurfaces in Toric Varieties}, hep-th/9506091,
Nucl. Phys. B456 (1995) 153.}
\ref\BV{N. Berline and M. Vergne, {\it Classes caracte\'eristiques
equivariantes, Formula de localisation en cohomologie
equivariante.} 1982.}
\ref\Bott{R. Bott, {\it A residue formula for holomorphic
vector fields}, J. Diff. Geom. 1 (1967) 311-330.}
\ref\Brion{M. Brion, {\it Equivariant cohomology and equivariant
intersection theory}, math.AG/9802063.}
\ref\CDGP{P. Candelas, X. de la Ossa, P. Green, and
L. Parkes, {\it A pair of Calabi-Yau manifolds as an
exactly soluble superconformal theory},
Nucl. Phys. B359 (1991) 21-74.}
\ref\CDFKM{P. Candelas, X. de la Ossa, A. Font, S. Katz
and D. Morrison, {\it Mirror symmetry for two parameter 
models I}, hep-th/9308083, Nucl. Phys. B416 (1994) 481-538.}
\ref\Cox{D. Cox, {\it
 The functor of a smooth toric variety}, alg-geom/9312001.} 
\ref\Delzant{T. Delzant, {\it Hamiltoniens p\'eriodiques
et image convex de l'application moment}, Bull. Soc.
Math. France 116 (1988) 315-339.}
\ref\EG{D. Edidin and W. Graham,{\it 
 Equivariant intersection theory}, alg-geom/9609018.}
\ref\Yau{S.T. Yau, ed.,{\it Essays on Mirror Manifolds I},
International Press, Hong Kong 1992.}
\ref\Fulton{W. Fulton, {\it Intersection Theory}, Springer-Verlag,
2nd Ed., 1998.}
\ref\FO{K. Fukaya and K. Ono, {\it
Arnold conjecture and Gromov-Witten invariants,}
preprint, 1996.}
\ref\GKZ{I. Gel'fand, M. Kapranov and A. Zelevinsky,
{\it Hypergeometric functions and toral manifolds},
Funct. Anal. Appl. 23 (1989) 94-106.}
\ref\GiventalI{A. Givental, {\it A mirror theorem for
toric complete intersections}, alg-geom/9701016.}
\ref\GKM{M. Goresky, R. Kottwitz and R. MacPherson,
{\it Equivariant cohomology, Koszul duality and the
localization theorem}, Invent. Math. 131 (1998) 25-83.}
\ref\GS{V. Guillemin and S. Sternberg, {\it Convexity properties
of the moment mapping}, Invent. Math. 67 (1982) 491-513.}
\ref\GSII{V. Guillemin and S. Sternberg, {\it Supersymmetry and
equivariant de Rham theory}, preprint.}
\ref\GZ{V. Guillemin and C. Zara, {\it Equivariant de Rham
Theory and Graphs}, math.DG/9808135.}
\ref\Hirzebruch{F. Hirzebruch, {\it Topological methods in
algebraic geometry}, Springer-Verlag, Berlin 1995, 3rd Ed.}
\ref\HKTY{S. Hosono, A. Klemm,
           S. Theisen and S.T. Yau, {\it Mirror symmetry, mirror map
           and applications to complete intersection Calabi-Yau
           spaces}, hep-th/9406055, Nucl. Phys. B433 (1995) 501-554.}
\ref\HLYII{S. Hosono, B.H. Lian and S.T. Yau, {\it GKZ-generalized
hypergeometric systems and mirror symmetry of Calabi-Yau hypersurfaces},
alg-geom/9511001, Commun. Math. Phys. 182 (1996) 535-578.}
\ref\HLY{S. Hosono, B.H. Lian and S.T. Yau, {\it Maximal
Degeneracy Points of GKZ Systems}, alg-geom/9603014,
Journ. Am. Math. Soc., Vol. 10, No. 2 (1997) 427-443.}
\ref\HSS{S. Hosono, M. Saito and J. Stienstra,
{\it On mirror conjecture for Schoen's Calabi-Yau 3-folds},
preprint 1997.}
\ref\Kr{A. Kresch,{\it Cycle groups for Artin Stacks,} math. AG/9810166.}
\ref\Kontsevich{M. Kontsevich,
{\it Enumeration of rational curves via torus actions.}
In: The Moduli Space of Curves, ed. by
R\. Dijkgraaf, C\. Faber, G\. van der Geer, Progress in Math\.
vol\. 129, Birkh\"auser, 1995, 335--368.}
\ref\LMW{ W. Lerche, P. Mayr, N. Warner,{\it Noncritical Strings, Del Pezzo 
Singularities and Seiberg-Witten Curves}, hep-th/9612085.}
\ref\LiTianII{ J. Li and G. Tian, {\it
Virtual moduli cycle and 
Gromov-Witten invariants of algebraic varieties},
J. of Amer. math. Soc. 11, no. 1, (1998) 119-174.}
\ref\LiTianIII{ J. Li and G. Tian, {\it
Virtual moduli cycle and Gromov-Witten invariants of general
symplectic manifolds.}
{\it Topics in symplectic $4$-manifolds 
(Irvine, CA, 1996), 47-83, First
Int. Press Lect. Ser. I}, Internat. Press,
Cambridge, MA, 1998.
}
\ref\LiTianIV{J. Li and G. Tian, {\it
Comparison of algebraic and symplectic GW-invariants},
To appear in Asian Journal of Mathematics, 1998.
}
\ref\LLY{B. Lian, K. Liu and S.T. Yau, {\it Mirror Principle I},
Asian J. Math. Vol. 1, No. 4 (1997) 729-763.}
\ref\LLYSurvey{B. Lian, K. Liu and S.T. Yau, {\it Mirror Principle, 
A Survey}, to appear.} 
\ref\MP{D. Morrison and R. Plesser, {\it
Summing the instantons: quantum cohomology and mirror symmetry
in toric varieties}, alg-geom/9412236.}
\ref\Oda{T. Oda, {\it Convex Bodies and Algebraic Geometry},
Series of Modern Surveys in Mathematics 15,
Springer Verlag (1985).}
\ref\Roan{S.-S. Roan, {\it The mirror of Calabi-Yau orbifolds},
Intern. J. Math. 2 (1991), no. 4, 439-455.}
\ref\Ruan{Y.B. Ruan, {\it Virtual neighborhoods 
and pseudo-holomorphic curves}, math.AG/9611021.}
\ref\RuanTian{Y.B. Ruan and G. Tian, {\it A mathematical
theory of quantum cohomology}, Journ. Diff. Geom. Vol. 42,
No. 2 (1995) 259-367.}
\ref\SiII{B. Siebert, {\it
Gromov-Witten invariants for general symplectic manifolds,}
alg-geom/9608005.}
\ref\Witten{E. Witten, {\it Phases of N=2 theories in two
dimension}, hep-th/9301042.}

\Title{}{Mirror Principle II} 

\centerline{
Bong H. Lian,$^{1}$\footnote{}{$^1$~~Department of Mathematics,
Brandeis University, Waltham, MA 02154.}
Kefeng Liu,$^{2}$\footnote{}{$^2$~~Department of Mathematics,
Stanford University, Stanford, CA 94305.}
and Shing-Tung Yau$^3$\footnote{}{$^3$~~Department of Mathematics,
Harvard University, Cambridge, MA 02138.} }

\vskip .2in

{\it Dedicated to Professor Michael Atiyah.}

\vskip .2in
Abstract. 
We generalize our theorems in {\it Mirror Principle I}
to a class of balloon manifolds. Many of the results
are proved for convex projective manifolds.
In a subsequent paper,
Mirror Principle III, 
we will extend the results to projective manifolds
without the convexity assumption.

\Date{}


\newsec{Introduction}

For the long history of mirror symmetry, consult \Yau.
For a brief description of more recent development, 
see the introduction in \LLY
\LLYSurvey. 
The present paper is a sequel to {\it Mirror Principle I} \LLY.
Here, we generalize all the results there to
a class of $T$-manifolds which we call balloon manifolds.
These results were announced in \LLYSurvey.

Let $X$ be a projective $n$-fold, and $d\in H^2(X,\Z)$.
Let $M_{0,k}(d,X)$ denote the moduli space of
$k$-pointed, genus 0, degree $d$, stable maps
$(C,f,x_1,..,x_k)$ on $X$
\Kontsevich. 
Note that our notation is without the bar.
By the work of \LiTianII
(cf. \BehrendFentachi\FO),
each nonempty  
$M_{0,k}(d,X)$ admits a cycle class $LT_{0,k}(d,X)$ in the Chow group
of degree $dim~X+\bra c_1(X),d\ket +n-3$. This
cycle plays the role of the fundamental class in
topology, hence
$LT_{0,k}(d,X)$ is
called the virtual fundamental class.

Let $V$ be a convex vector bundle on $X$.
(ie. $H^1(\P^1,f^*V)=0$ for every holomorphic
map $f:\P^1\ra X$.) Then $V$ induces on each
$M_{0,k}(d,X)$ a vector bundle $V_d$,
with fiber at
$(C,f,x_1,..,x_k)$ given by the section space $H^0(C,f^*V)$.
Let $b$ be any multiplicative characteristic class
\Hirzebruch.
(ie. if $0\ra E'\ra E\ra E''\ra 0$ is an exact sequence
of vector bundles,
then $b(E)=b(E')b(E'')$.) 
The problem we study here is to compute
the characteristic numbers
$$K_d:=\int_{LT_{0,0}(d,X)} b(V_d)$$
and their generating function:
$$\Phi(t):=\sum K_d~ e^{d\cdot t}.$$
There is a similar and equally important problem if one
starts from a concave vector bundle $V$ \LLY.
(ie. $H^0(\P^1,f^*V)=0$ for every holomorphic
map $f:\P^1\ra X$.) More generally, $V$ can
be a direct sum of a convex and a concave bundle.
Important progress made on these problems has come from
mirror symmetry. All of it
seems to point toward the following
general phenomenon \CDGP, which we call
{\it the Mirror Principle}. Roughly,
it says that the function $\Phi(t)$ can be computed
by a change of variables in terms of certain explicit
special functions, loosely called generalized
hypergeometric functions.

When $X$ is a toric manifold with $c_1(X)\geq0$,
$b$ is the Euler class, and $V$ is a sum
of line bundles, there is a general formula
derived in \HLY~ from mirror symmetry.
This formula was later studied in \GiventalI~
based on a series of
axioms.

{\bf Acknowledgements.} We thank B. Cui, J. Li, T.-J. Li, and G. Tian,
for numerous helpful discussions with us during the course
of this project. B.H.L.'s research
is supported by NSF grant DMS-9619884.
K.L.'s research is supported by NSF grant
DMS-9803234 and the Terman fellowship.
S.T.Y.'s research is supported by DOE grant
DE-FG02-88ER25065 and NSF grant DMS-9803347.

\subsec{Main Ideas}

We now sketch our main ideas for computing the classes $b(V_d)$.

{\it Step 1. Localization on the linear sigma model.} 
Consider the moduli spaces 
$M_d(X):=M_{0,0}((1,d),\P^1\times X)$. The projection
$\P^1\times X\ra X$ induces a map $\pi:M_d(X)\ra M_{0,0}(d,X)$.
Moreover, the standard action of $S^1$ on $\P^1$ induces
an $S^1$ action on $M_d(X)$. We first study a 
slightly different problem. Namely
consider the classes $\pi^*b(V_d)$ on $M_d(X)$,
instead of $b(V_d)$ on $M_{0,0}(d,X)$.
First, there
is a canonical way to
embed fiber products (see below)
$$F_r=M_{0,1}(r,X)\times_X M_{0,1}(d-r,X)$$
each as an $S^1$ fixed point component into $M_d(X)$.
Let $i_r:F_r\ra M_d(X)$ be the inclusion map.
Second, there is an evaluation map
$e:F_r\ra X$ for each $r$.
Third, suppose that
there is a projective manifold $W_d$ with $S^1$ action,
that there is an equivariant map $\varphi:M_d(X)\ra W_d$, and
embeddings $j_r:X\ra W_d$, such that the diagram
$$\matrix{
F_r & {\br i_r\over \longrightarrow} & M_d(X)\cr
e\downarrow &  & \downarrow \varphi\cr
X & {\br j_r\over \longrightarrow} & W_d}
$$
commutes.
Let $\alpha$ denotes the weight of
the standard $S^1$ action on $\P^1$.
Then applying the localization formula \AB, this diagram allows us
to recast our problem to one of studying
the $S^1$-equivariant classes
$$Q_d:=\varphi_!\pi^*b(V_d)$$
defined on $W_d$.
Moreover we can expand the class 
$$A_d:={j^*_0Q_d\over e_{S^1}(X_0/W_d)}$$
 on $X$
in powers of $\alpha^{-1}$,
and find that it is of order $\alpha^{-2}$.

The spaces $W_d$ in the commutative diagram above
are called the linear sigma model of $X$.
They have been introduced in \MP~following \Witten~
when $X$ is a toric manifold,

{\it Step 2. Gluing identity.}
Consider the vector bundle $\cU_d:=\pi^* V_d\ra M_d(X)$,
restricted to the fixed point components $F_r$.
A point in $(C,f)$ in $F_r$ is
a pair $(C_1,f_1,x_1)\times(C_2,f_2,x_2)$ of
1-pointed stable maps glued together at
the marked points, ie. $f_1(x_1)=f_2(x_2)$.
From this, we get an exact sequence of bundles
on $F_r$:
$$0\rightarrow i_r^*\cU_d\rightarrow U_r'\oplus U_{d-r}'\rightarrow
e^*V\rightarrow 0.$$
Here $i_r^*\cU_d$ is the restriction to $F_r$ of the bundle
 $\cU_d\ra M_d(X)$. And $U_r'$ is the pullback of the bundle
$U_r\ra M_{0,1}(d,X)$ induced by $V$, and similarly for $U_{d-r}'$.
Taking the multiplicative characteristic class $b$, we get
the identity on $F_r$:
$$e^*b(V)b(i^*_r\cU_d)=b(U_r')b(U_{d-r}').$$ 
This is what we call
the {\it gluing identity}.
This may be translated to a similar quadratic identity, via
Step 1, for $Q_d$ in the equivariant cohomology 
groups $H^*_{S^1}(W_d)$.
The new identity is called the Euler data identity.

{\it Step 3. Linking theorem.}
The construction above is functorial, so that if $X$
comes equipped with a torus $T$ action, then the entire
construction becomes $G=S^1\times T$ equivariant and
not just $S^1$ equivariant. In particular,
the Euler data identity is an identity of
$G$-equivariant classes on $W_d$. 
Our problem is  to first
compute the $G$-equivariant classes $Q_d$ on $W_d$
satisfying the Euler data identity, and with the property that
$A_d\sim\alpha^{-2}$.
Note that  the restrictions $Q_d|_p$ to the $T$ fixed points
$p$ in $X_0\subset W_d$ are polynomial functions on
the Lie algebra of $G$.
Suppose that $X$ is a balloon manifold. 
Then it can be shown that (with a nondegeneracy assumption
on $e_G(X_0/W_d)$) the classes $Q_d$ are uniquely determined
by the values of the $Q_d|_p$, when $\alpha$ is 
some scalar multiple of a weight on the tangent space
$T_pX$. These values of $Q_d|_p$ can be
computed explicitly by exploiting the
structure of a balloon manifold.

Once these values are known, it is often easy to manufacture
explicit $G$-equivariant classes $\tilde Q_d$ with 
the restrictions $\tilde Q_d|_p$ having the above same values,
and satisfying the Euler data identity.
In this case, we say that the data $\tilde Q_d$
are linked to the data $Q_d$.
By a suitable change of variables, one
can also arrange that
${j^*_0\tilde Q_d\over e_{S^1}(X_0/W_d)}\sim \alpha^{-2}$.
By the preceding discussion, we get $Q_d=\tilde Q_d$.

{\it Step 4. Computing $\Phi(t)$.}
Once the classes $Q_d=\varphi_!\pi^* b(V_d)$ are determined,
we can unwind the many maps used in Step 1. 
The preceding computations can be done simply in
the form of power series. This finally computes
the generating function $\Phi(t)$. 

\bs
The answer for $\Phi(t)$ is given in the form of Conjecture 9.1.
In this paper, for clarity, we restrict ourselves to the
case when the tangent bundle of $X$ is convex. We prove that
Conjecture 9.1 holds whenever $X$ is a balloon manifold
having a linear sigma model $W_d$ such that $e_G(X_0/W_d)$ 
satisfies a nondegeneracy condition.

In the nonconvex case, we must replace $M_{0,k}(d,X)$
by Li-Tian's virtual fundamental cycle \LiTianII~
for the purpose of localization and integration. 
The sequel, Mirror Principle III, to this paper
will be devoted entirely to dealing with the added
technicality arising from this replacement.
All the results in this paper will 
generalize with only slight modifications
as a result of 
this replacement, but with no change to
the overall conceptual framework.

By the equivalence,
established in \LiTianIII,
of symplectic GW theory and algebraic GW theory for
projective manifolds,
we also expect that the results in this paper can be readily
generalized to the symplectic case \RuanTian\SiII.

\newsec{Set-up}

\subsec{Equivariant localization} 

We first discuss some basic facts about localization. 
The key technique of our proof is the equivariant
localization formula, due to  Atiyah-Bott
\Atiyah\Bott\AB, and Berline-Vergne \BV. 
For an orbifold version of the localization formula, see \Kr. The
spirit of the localization we'll use
is closer to the Bott residue formula. 
We first explain
this formula.

Let $X$ and $Y$ be two spaces, by which we mean compact manifolds or orbifolds, with a torus
$T$-action. When an orbifold is involved, the integral and localization formulas should be taken in the orbifold sense. Let $\{ F\} $ be the components of the fixed point set. Let
$H^*_T(\cdot)$ denote the equivariant cohomology group with complex coefficient,
and $i_F:F\ra X$ the inclusion map. We say
that equivariant localization
holds on $X$, if the two maps
$$i_F^*: \ \ H^*_T(X)\rightarrow H^*_T(F),\ \ {i_F}_!:\ \ 
H^*_T(F)\rightarrow \ H^*_T(X).$$
which are respectively the pull-back, and the Gysin map, are such
that the following formulas holds:
given any equivariant cohomology class $\omega$ on $X$, we have
$$\omega=\sum_F{i_F}_!\left({i_F^*\omega\over e_T(F/X)}\right).$$
This formula is equivalent to the integral version of the localization
formula
$$\int_X \omega
=\sum_F\int_F{i_F^*\omega\over e_T(F/X)}.$$

An important fact about equivariant theory
is that, if $V$ is an equivariant vector bundle
on an orbifold $X$, then any characteristic class of $V$ has an
equivariant extension.  
Let $T=S^1$
for simplicity.
If $c_{2k}$ is a
characteristic class of degree $2k$, then its equivariant
extension can be represented by the form
$$c=c_{2k}+c_{2k-2}+\cdots +c_0$$ 
in the equivariant cohomology of $X$.
 
Here is a way to calculate the terms in the localization
formula. Assume that
each fixed point component $F$ is smooth. 
If $c_{2k}$ is a Chern class, then by using splitting
principle, it can be  expressed
as a symmetric function of the Chern roots:
 $P(x_1, \cdots, x_l)$ where
$l=rank~ V$. When $V$ can be
decomposed into a direct sum of line bundles on $F$:
$$V|_F=L_1\oplus \cdots \oplus L_l$$ 
with the $T$-action on $L_j$ given by,
say the character
$e^{2\pi \sqrt{-1} n_jt}$, then the restriction of its equivariant
counterpart $c$
to $F$ is
$$i_F^*c=P(c_1(L_1)+n_1t, \cdots, c_1(L_l)+n_lt).$$
 
The computation of the equivariant Euler class of $F$ in $X$ is similar.
When the restriction of $TX$ to $F$ has a decompostion into line
bundles
$$TX|_F=E_1\oplus \cdots \oplus E_n$$ 
where $T$ acts on $E_j$ by the character
$e^{2\pi \sqrt{-1}m_jt}$, then
$$e_T(F/X)=\prod_j(c_1(E_j)+m_jt).$$
In the above, the $n_j$'s and $m_j$'s are integers for if $X$ is a manifold, and are ratonal numbers if $X$ is an orbifold.

\subsec{Functorial localization formula}

In this subsection, we derive two
formulas which are often used in our work.
Let  $X,Y$ be two $T$-spaces, and
$$f:\ \ X \rightarrow Y$$
be an equivariant map. Let $E$ be a fixed point component in $Y$, and
$F:=f^{-1}(E)$ be a fixed component in $X$. Let $g$ be the restriction
of $f$ to $F$, and $j_E:E\ra Y$, $i_F:F\ra X$ be the inclusion maps.
Thus we have the commutative diagram:
$$\matrix{
F & {\br i_F\over \longrightarrow} & X\cr
g\downarrow &  & \downarrow f\cr
E & {\br j_E\over \longrightarrow} & Y.}
$$

\lemma{ Given any class $\omega \in H^*_T(X)$, we have the
equality on $E$:
$${j_E^*f_!(\omega)\over e_T(E/Y)}
=g_!\left({i_F^*\omega\over e_T(F/X)}\right)
$$
}
\proof Let us consider localization of $\omega f^*{j_E}_!(1)$ on $X$,

$$\omega f^*{j_E}_!(1)={i_F}_!\left(
{i_F^*(\omega
f^*{j_E}_!(1))\over e_T(F/X)}\right).$$ 
Note the contributions from fixed
components other than $F$ vanish. Applying
 the push-forward $f_!$ to both
sides, we get

$$f_!(\omega){j_E}_!(1)=f_!{i_F}_!
\left({i_F^*(\omega
f^*{j_E}_!(1))\over e_T(F/X)}\right).$$
Now $f\circ i_F=j_E\circ g$ which, implies
$$f_!{i_F}_!={j_E}_!g_!,\ \ i_F^* f^*=g^*j_E^* .$$ 
Thus we get
$$f_!(\omega){j_E}_!(1)={j_E}_!g_!
\left({i_F^*(\omega)~g^*e_T(E/Y)\over e_T(F/X)}\right).$$
Applying $j_E^*$ to both sides, we then arrive at

$$j_E^*f_!(\omega)~e_T(E/Y)=e_T(E/Y)~
g_!\left({i_F^*(\omega)~g^*e_T(E/Y)\over e_T(F/X)}\right)=
e_T(E/Y)^2~g_!\left({i_F^*(\omega)\over e_T(F/X)}\right).$$
Since $e_T(E/Y)$ is invertible, our assertion follows.
$\Box$

The same argument applies to the case when $E$ and $F$ are $T$-invariant subspaces. A slightly different argument for the proof of the above lemma 
will be given in our subsequent paper. We will also
need the following formula, 
which is actually a special case of \Fulton, Theorem 6.2. Here we include a proof for the convenience of the reader. The spaces involved are $T$-spaces, that is, $T$-manifolds or orbifolds.

\lemma{Suppose we have a $T$-equivariant commutative diagram
$$\matrix{
V & {\br i\over \longrightarrow} & W\cr
g\downarrow &  & \downarrow f\cr
Z & {\br j\over \longrightarrow} & Y}
$$
such that $f^* j_!(1)=i_!(1)$. Then for any class
$\omega$ on $H_T^*(W)$, we have the following equality
on $Z$:
$$j^*f_!(\omega)=g_! i^*(\omega).$$
}
\proof
By assumption, we have
$$\omega~f^*j_!(1)=\omega~i_!(1)=i_!i^*(\omega).$$
Applying $j^*f_!$ to both sides, on the one hand we get
$$
j^*f_!(\omega~f^*j_!(1))
=j^*(f_!(\omega)~j_!(1))
=j^*f_!(\omega)~e_T(Z/Y).
$$
On the other hand, we get
$$
j^*f_!(i_!i^*(\omega))
=j^*j_!g_!i^*(\omega)
=g_!i^*(\omega)~e_T(Z/Y).
$$
Thus our assertion follows.  $\Box$

The case we will use in this paper is when $Z$ and $V=f^{-1}(Z)$ are both $T$ invariant submanifolds of same codimension, in which the condition in the Lemma clearly holds. 

\subsec{Balloon manifolds}

By a {\it balloon manifold}, we mean a complex projective
$T$-manifold $X$ with the following properties. There are
only finite number of $T$-fixed points. At each fixed point $p$,
the $T$-weights on the isotropic representation $T_pX$ are
pairwise linearly independent. This class of
manifolds were introduced by Goresky-Kottwitz-MacPherson
 \GKM. (We refer the reader to \GSII~ for an
excellent exposition.)
Throughout this paper,
we assume that {\it $X$ is convex}, ie. 
$H^1(\P^1,f^*TX)=0$ for any holomorphic map $f:\P^1\ra X$.

One important property of 
a balloon $n$-fold is that at each fixed point $p$,
there are exactly $n$ balloons, ie. $T$-invariant $\P^1$,
each balloon connecting $p$ to one other fixed point $q$.
The induced action on each balloon is
the standard rotation with two fixed points $p$ and $q$. 
(see \GKM \GZ). 
We {\it denote by $pq$ the balloon connecting the fixed points $p,q$.}
Toric manifolds, complex $C$-spaces and spherical manifolds are examples
of balloon manifolds.

We fix a $T$ equivariant embedding of $X$ into the product of projective
spaces 
$$\P(n):=\P^{n_1}\times \cdots \times \P^{n_m}$$
 such that the pull-backs of the hyperplane classes $H=(H_1, \cdots,
H_m)$ generate $H^2(X, {\bf Q})$. We use the same notations for the
corresponding equivariant classes of the $H$'s, and their restrictions
to $X$. For $\omega\in H^2(X)$ and $d\in H_2(X)$, we
denote their pairing by $\bra\omega, d\ket$. 

For convenience, we introduce the following notations:
$$\eqalign{
H&=(H_1,..,H_m)\cr
H\cdot\zeta&=H_\zeta=H_1\zeta_1+\cdots+H_m\zeta_m\cr
H(p)&=(H_1(p),..,H_m(p))\cr
H_\zeta(p)&=H_1(p)\zeta_1+\cdots+H_m(p)\zeta_m.
}$$
Here $\zeta=(\zeta_1,..,\zeta_m)$ are formal variables.
We denote by $K^\vee\subset H_2(X)$ the 
set of points in $H_2(X,\Z)_{free}$ 
in the dual of the closure of the K\"ahler cone of $X$.
Since $K^\vee$
is a semigroup in $H_2(X)$, it defines
a partial ordering $\succ$ on the lattice $H_2(X,\Z)_{free}$.
That is, $d\succeq r$ iff $d-r\in K^\vee$.
Let $\{H_j^\vee\}$ be the basis dual to the $\{H_j\}$ in $H_2(X)$.
If $d\succeq r$ for two classes $d,\, r\in
H_2(X)$, then $d-r= d_1H^\vee_1+\cdots +d_mH^\vee_m$ for nonnegative
integers $d_1, \cdots, d_m$.

We also consider a balloon manifold as a symplectic manifold with
a symplectic structure
given by $\omega=H_\zeta$ for some generic $\zeta$. 
By the
convexity
theorem of Atiyah \Atiyah~ and Guillemin-Sternberg \GS, 
the image of the moment map $\mu_\zeta$ in the dual Lie algebra ${\cal T}^*$ is a convex polytope,
known as
the moment polytope.
When $X$ is a toric manifold, the moment polytope
is known as a Delzant polytope \Delzant. In this case, it is
well-known that the normal fan of this polytope
is the defining fan of $X$. 

We say $X$ a multiplicity-free manifold, if for each point $p$ in ${\cal T}^*$, the inverse image $\mu_\zeta^{-1}(p)$ is connected.

\lemma{ Let $X$ be a multiplicity-free balloon
manifold, then $H(p)\neq
H(q)$ for any two distinct fixed points $p$ and $q$ in $X$. }
\proof  Let $\mu_\zeta$ denote the moment map of the $T$-action on $X$
with
respect to the symplectic form $H_\zeta =\bra H, \zeta\ket$ 
for a generic
choice of $\zeta\in \C^k$. Then the image of $\mu_\zeta:\
X\rightarrow \cT^*$ is a convex polytope whose vertices are given the
images
of the fixed points $\{ p\}$. The weights of
$H_\zeta $ at the fixed points, up to an over all translation, are the
same as $\{
\mu_\zeta (p)\}$ which are all different. Since $X$ is
multiplicity-free, the inverse image of each vertex contains only one
fixed point in $X$. Since $\zeta$ is generic, this implies
that the $H(p)$'s are distinct at different fixed point.  $\Box$

We shall assume throughout this paper that 
{\it $H(p)\neq H(q)$ for all distinct fixed points $p,q$ in $X$.}
Equivalently, 
if $c(p)=c(q)$ for all $c\in H^2_T(X)$, then $p=q$.
This condition is also equivalent to
the statement that the moment map with respect to $\omega=H_\zeta$ and the $T$
action is injective to the set of vertices of the moment polytope, when
restricted to the fixed points $X^T$. By the above lemma we know
that  toric manifolds and compact homogeneous manifolds all satisfy this
condition.

When $X$ is a toric $n$-fold, we have $N=m+n$
$T$-invariant divisors in $X$. Let
$D_a=c_1(L_a)$, $a=1,..,N$, be the equivariant
first Chern classes of the corresponding 
equivariant line bundles. 
These $T$ divisors correspond 1-1 with
the one-cones of the defining fan of $X$ \Oda.
Moreover the fixed points correspond 1-1 with the $n$-cones.
Labelling the $n$-cones by $p\in X^T$, we have
a balloon \Oda~ $pq$ in $X$ iff the $n$-cones
$p,q$ intersect in a codimension one subcone. Since $X$
is smooth, hence the $n$-cones are regular, there are
exactly $n$ balloons $pq$ for each fixed $q$.
One can give a dual description of all these
by using the Delzant polytope.

Returning to the general case,
suppose that $X$ is a balloon manifold,
and that we have equivariant classes $\{D_a\}$
in $H^2_T(X)$ with the following
property.
At every fixed point $p$, $D_a(p)$ is either zero or
it is a weight
on $T_pX$.
Let $pq$ be a balloon in $X$.
The induced $T$-action on $pq$ is
the standard rotation with fixed points $p,\, q$. By applying
the localization formula on $pq\simeq \P^1$ and the integral 
$\bra c,[pq]\ket$,
we have
$$c(q)=c(p)+\bra c,[pq]\ket\lambda$$
for all $c\in H^*_T(X)$, 
where $\lambda$ is the weight on the tangent line $T_q(pq)$.
Let $\lambda=D_a(q)$. Specializing to $c=D_a$,
we get
$D_a(q)=D_a(p)+\bra D_a,[pq]\ket D_a(q)$. This shows that
$\bra D_a,[pq]\ket\neq0$. For otherwise we would have
$D_a(q)=D_a(p)\neq0$, and this would mean that $D_a(p)$
is a weight on $T_p(po)$ for some edge $po$
running {\it in the direction} of 
$D_a(q)=D_a(p)$ from $p$ to $o$.
So  we had 
three vertices lying joined in a line from $q$ to $p$ to $o$ 
in the moment graph. This would mean that there is a pair of
linearly dependent weights on the tangent space $T_pX$,
which can't happen in a balloon manifold.
A similar argument shows that
$\bra D_a,[pq]\ket=1$.

\lemma{Let $\omega=H_\zeta$ and $p,q\in X^T$, $r\succ0$ and
$\lambda$ be a weight on $T_qX$.
If $\omega(q)=\omega(p)+\bra \omega,r\ket\lambda$ for generic $\zeta$,
then $p,q$ are joined by a balloon, $r=[pq]$, and
$\lambda$ is the weight on the tangent line $T_q(pq)$.}
\thmlab\BalloonLemma
\proof
It suffices to prove that $p,q$ are joined by a balloon.
The last two conclusions then follow immediately. 
Then under the corresponding moment map $\mu$, $p,q$ are
mapped to $\omega(p), \omega(q)$ (up to an overall affine
transformation), which are distinct
because
\eqn\dumb{\omega(q)-\omega(p)=\bra \omega,r\ket\lambda\neq0.}
Since $\lambda$ is a weight on $T_qX$, there is an edge emanating from 
the point $\omega(q)$ in the direction of $\lambda$, ending at some
other vertex $\omega(o)$, where $qo$ is a balloon in $X$.
If $\omega(p)\neq \omega(o)$,
we would have three distinct vertices
of the moment polytope lying on a single line.
Thus $\omega(p)=\omega(o)$, which implies that $p=o$.  $\Box$

\lemma{
The zero class $\omega=0$ is the only class in $H^*_T(X)$ with
the property that
$$\int_X \omega~ e^{H_\zeta}=0$$
for all generic $\zeta \in \C$.
}\thmlab\NondegeneracyLemma
\proof
Suppose
$$\int_{X}\omega~ e^{H_\zeta}=0.$$
By localization, we have
$$\sum_{p\in X^T}{\omega(p)
\over e_T(p/X)} e^{H_\zeta(p)}=0.$$
But since the vectors $H(p)$ are all distinct,
those exponential functions in $\zeta$ are linearly 
independent over the field $\Q(\cT^*)$,
implying that
$\omega(p)=0$ for all $p$. Thus  
$\omega=0$. 
$\Box$

\subsec{Sigma models}
\seclab\SigmaModels

Let $X$ be balloon manifold
with a fixed $T$-equivariant embedding
$X\ra\P(n)$, as discussed above.
We write
$$M_d(X):=M_{0,0}((1,d),
\P^1\times X).$$
Since $X$ is assumed to be
convex, $M_d(X)$ is an orbifold. The standard $S^1$ action 
on $\P^1$ together with the $T$ action on $X$ induce
a $G=S^1\times T$ action on $M_d(X)$.

Here is a description of some $S^1$ fixed point components $F_r$,
labelled by $0\preceq r\preceq d$, inside of $M_d(X)$.
Let $F_r$ be the fiber product
$$F_r:= M_{0,1}(r, X)\times_X M_{0,1}(d-r, X)$$
More precisely, consider the map
$$ev_r\times  ev_{d-r}: 
M_{0,1}(r, X)\times M_{0,1}(d-r, X)\rightarrow X\times X$$ 
given by evaluations at the marked points; and
$$\Delta:\ X\rightarrow X\times X$$ 
the diagonal map. Then 
$$F_r=(ev_r\times ev_{d-r})^{-1}\Delta(X).$$
Note that $F_d=M_{0,1}(d,X)$ by convention.
The set $F_r$ can be identified with an $S^1$ fixed point
component of $M_d(X)$ as follows. 
Consider the case $r\neq0,d$ first. Given a point
$(C_1,f_1,x_1)\times (C_2,f_2,x_2)$ in $F_r$, we
get a new curve $C$ by gluing
$C_1,C_2$ to $\P^1$ with
$x_1,x_2$ glued to $0,\infty\in\P^1$ respectively.
The new curve $C$ is mapped into $\P^1\times X$ as follows.
Map $\P^1\subset C$ identically onto $\P^1$, and
collapse $C_1,C_2$ to $0,\infty$ respectively; 
then map $C_1,C_2$ into $X$ with $f_1,f_2$ respectively,
and collapse the $\P^1$ to $f(x_1)=f(x_2)$. 
This defines a point
$(C,f)$ in $M_d(X)$. 
For $r=0$, we glue $(C_1,f_1,x_1)$ to
$\P^1$ at $x_1$ and $0$. For $r=d$, we glue
$(C_2,f_2,x_2)$ to $\P^1$ at $x_2$ and $\infty$.
We will identify $F_r$
as a subset of $M_d(X)$ as above, and let
$$i_r:F_r\ra M_d(X)$$
denotes  the inclusion map. 
Clearly, we also have an evaluation map 
$$e_r:F_r\ra X$$
which sends a pair in $F_r$ to the value
at the marked point. In the following, we will simply write $e_r$ as $e$ without causing any coonfusion.

We call a compact manifold or orbifold $W_d$ with $G=S^1\times T$ action a {\it linear sigma model} of degree $d$ for $X$, if the following conditions are satisfied:

\item{1.} The $S^1$ action on $W_d$ has fixed point components given by $X_r$, labelled by $0\preceq r\preceq d$, and each $X_r$ is $T$-equivariantly isomorphic to $X$.

\item{2.} There is a $G$-equivariant birational map $\varphi$ from $M_d(X)$ to $W_d$, such that $\varphi|_{F_r}=e$, and $\varphi^{-1}(X_r)=F_r$.

\item{3.} All equivariant cohomology classes in $H^2_G(W_d)$ are lifted from $H_T^2(X)$, and the lift $\hat{D}\in H_G^2(W_d)$ of $D\in H_T^2(X)$ restricts to $D+\bra D, r\ket\alpha$ on $X_r$. 

\item{4.} The $G$-equivariant Euler class of the normal bundle of $X_0$ in $W_d$ has the form
$$e_G(X_0/W_d)= \prod_a \prod_{m_a}(D_a-m_a \alpha)$$
where the $m_a$'s are positive integers and the $D_a$'s are classes in $H^2_T(X)$, such that at a given $T$ fixed point $p$ in $X$,
the nonzero $D_a(p)$'s are multiples of distinct weights of $T_pX$.

\bs
Here a birational map, in algebraic geometry language, is a regular morphism which is an isomorphism when restricted to a Zariski open set in $M_d(X)$.

Note  $W_d$ need not be unique. We identify $X_r$ with $X$ by assumption
 1, and denote by
$$j_r:X_r\ra W_d$$
the inclusion map. 

We call a balloon manifold $X$ {\it admissible}
if it has a linear sigma model $W_d$ for each $d$, 
and that $H_\zeta(p)\neq H_\zeta(q)$ 
for any two distinct fixed points $p,q$ in $X$.
The main result in this paper is to show that the 
mirror principle holds for any admissible balloon manifold.

\remark{Condition 4 is actually assuming more than what we need.
This condition can be replaced by the following
weaker, but more technical condition.
For each fixed point $p$ and for any $d$, as a function of $\alpha$,
$e_G(X_0/W_d)|_p$ has possible zero only  at either 0
or a multiple of a weight $\lambda$ on $T_pX$. In addition if $[pq]$
is a balloon and $d=\delta[pq]$, then $\lambda/\delta$
is at worse a simple zero.
For example, the following form would meet this criterion:
$$
e_G(X_0/W_d)=
{\prod_a\prod_{m_a}(D_a-m_a\alpha) \over
\prod_b\prod_{n_b}(D_b-n_b\alpha)}
$$
where the $m_a,n_b$ are nonzero scalars.
}

{\it Example 1}: Projective space $\P^n$ with $W_d=\P^{(n+1)d+n}$ is 
admissible. 
The existence of $\varphi$ was proved in \LLY, which is the so-called Li-Tian map. The lifted hyperplane class $\kappa$ has the required property that 

$$j_r^*\kappa =H+\bra H, r\ket\alpha= H+r\alpha.$$ The equivariant Euler class

$$e_G(\P^n_0/W_d)=\prod_{i=0}^n\prod_{m=1}^d (H-\lambda_i-m\alpha)$$where $\lambda_i$'s denote the weights of the torus $T$ action on $\P^n$. Clearly the equivariant classes $\{ H-\lambda_i\}$ has the required property.

{\it Example 2}: More generally for $\P(n)$, we can take $N_d(\P(n))$ to be $W_d$. In fact the $S^1$ fixed point components on $N_{k,l}$
are exactly $k+1$ copies $\P^l_r$, $r=0,..,k$, of
$\P^l$. Each $\P^l_r$ consists of $l+1$ tuples
of monomials, each being a scalar multiple of $w_0^r w_1^{k-r}$.
Similarly the $S^1$ fixed point components on $N_d(\P(n))$
are copies $\P(n)_r$, $0\preceq r\preceq d$,
of $\P(n)$. All equivariant
cohomology classes in $H^2_G(N_d(\P(n))$ are lifted from $H^2_T(\P^n)$
(cf \LLY). Let $\kappa_i$ be the lift of the hyperplane class $H_i$, 
of the $i$th factor $\P^{n_i}$. Then
$$j_r^*\kappa_i=H_i+\bra H_i,r\ket\alpha$$
where $j_r$ denotes the inclusion of 
$\P(n)_r$ in $N_d(\P(n))$. By using the formula in \LLY, it is easy to show the equivariant Euler class $e_G(\P(n)_0/N_d(\P(n)))$, which is a product of $e_G$'s in last example, has the required property.

{\it Example 3:} Let $N_{k,l}$ be the space of $l+1$ tuples $[f_0,..,f_l]$
of degree $k$ polynomials $f_i(w_0,w_1)$, modulo scalar.
Thus $N_{k,l}\cong\P^{(l+1)k+l}$. It is called
the linear sigma model for $\P^l$. (See \LLY.)
Let
$$N_d(\P(n)):=N_{d_1,n_1}\times\cdots \times N_{d_m,n_m}.$$
Recall that we have a collapsing map $\varphi:M_k(\P^l)\ra 
N_{k,l}$, which is $G:=S^1\times T$ equivariant. By taking composite with the projection from $M_d(\P(n))$ to each $M_{d_j}(\P^{n_j})$, we obtain a $G$-equivariant map
$$M_d(\P(n))\ra N_d(\P(n))$$
which we also denote by $\varphi$. 
Note that $M_d(X)$ can be viewed as a cycle in $M_d(\P(n))$. 
We denote 
the image cycle $\varphi_!(M_d(X))$ in $N_d(\P(n))$ 
$\varphi$ by $ N_d(X)$.

If $N_d(X)$ is a manifold or an orbifold, then Properties 1-3 are automatically satisfied, if furthermore $e_G(X_0/W_d)$ has property 4, then we can simply take $W_d=N_d(X)$ as the linear sigma model.

{\it Example 4:} Convex toric varieties. In this case $W_d$ is a toric
$n$-fold, as introduced by Witten \Witten
and used first by Morrison-Plesser \MP to
study quantum cohomology. Recall that a toric manifold
$X$ can be realized as the GIT quotient $\C^N//T^m_{\C}$ 
where $T^m_\C$ is a $m$-dimensional complex torus
acting on $\C^N$. Here $m=rank~H^2(X,\Z)$, $N=n+m$.
Let $[z_1, \cdots,
z_N]$ denote the coordinates on $\C^N$.
Then each $z_j$ can be viewed as a section of a line bundle $L_j$
on $X$ \Cox\MP. Modulo the induced action by $T^m_C$ 
from $\C^N$, a map from $\P^1$ into $X$ is uniquely
represented by an $N$-tuple of polynomials 
$$[f_1(w_0, w_1), \cdots, f_N(w_0, w_1)]$$ 
where $f_j$ is a section of
the line bundle $O(l_j)$ over $\P^1$ with $l_j=\bra c_1(L_j), d\ket$. 
Let $\C^N(d)$ be the vector space of $N$-tuple
of polynomials of degree
$(d_1, \cdots, d_N)$ as above. Then as described in \MP, $W_d$ is the
GIT quotient by the induced action of $T^m_C$ on it: 
$$W_d=\C^N(d)//T^m_C.$$  

Let $M_d^o(X)$ denote the set of points $(f, C)$ in $M_d(X)$ such that $C\simeq \P^1$. We call $M_d^o(X)$ the smooth part of $M_d(X)$. We can define a map $\varphi_o$ from $M_d^o(X)$ to $W_d$ in the following way: each $(f, C)$ gives a map from $\P^1$ to $X$, and modulo the induced $T^m_C$ action, uniquely determines $N$-tuple of polynomials as above, therefore gives a point in $W_d$, which we define to be the image of $(f, C)$ under $\varphi_o$. 
This is clearly a canonical identification.

It is not difficult to see that the $S^1$-fixed components in $W_d$ can
be described as GIT quotient,
$$X_r\simeq \{[a_1 w_0^{\bra D_0, r\ket}w_1^{\bra D_0,d- r\ket},
\cdots, a_Nw_0^{\bra D_N,
r\ket }w_1^{\bra D_N,d- r\ket }]|a\in\C^N\}//T^m_C.$$ 
The equivariant Euler class of its
normal bundle in $W_d$ is
$$e_G(X_r/W_d)=\prod_{a=1}^N
\prod_{k=0, k\neq \bra D_a, r\ket }^{\bra D_a,
d\ket}(D_a+\bra D_a, r\ket\alpha-k\alpha).$$
Here $D_a=c_1(L_a)$ is the equivariant 
first Chern class of the line bundle $L_a$ corresponding to the
$a$th component in the coordinates of $X$. 
The lift $\hat{D}_a$ of $D_a$ to $W_d$ clearly has the property 

$$j_r^*\hat{D}_a=D_a+\bra D_a, r\ket\alpha.$$
As pointed out in \MP, the cohomology of $W_d$ are generated by
the $\hat D_a$.
Thus $W_d$ has properties 1, 3, and 4.
In the next subsection, we establish that the $\varphi_o$
extends to a regular $G$-equivariant
map $\varphi$ from $M_d(X)$ to $W_d$, with
property 2.
So, for a toric $X$ we can take its linear
sigma model to be $W_d$ as constructed above. 
It follows that $X$ is an admissible balloon manifold.

{\it Example 5:} Our method works well even for certain singular
manifolds. We take weighted projective space as example to illustrate
the ideas. Let $\P^n_a$ with $a=(a_0, \cdots, a_n)$ be a weighted
projective space. Let $[z_0, \cdots, z_n]$ be the coordinates for
$\P^n_a$, then $z_j$ can be considered as a section of the line bundle
$O(a_j)$. Then the linear sigma model for this weighted projective space
is the induced weighted quotient by $S^1$ on the space of $n+1$ tuple of
polynomials $[f_0(w_0, w_1), \cdots, f_1(w_0, w_1)]$ where $f_j$ is a
section of the line bundle $O(da_j)$ on $\P^1$. 

It is known that $\P^n_a$ is equivalent to $\P^n/Z_a$ where $Z_a$ is a
finite group. The space $M_d(\P^n_a)$ for $\P^n_a$ is
equal to $M_d(\P^n)/Z_a$. On the other hand we can also take 
$N_d(\P^n)/Z_a$ as the linear sigma model $W_d$. In this case $W_d$ is
an orbifold, a weighted projective space. Since the action of $Z_a$ commutes with the action of torus $T$, we see that the induced collapsing map 

$$\varphi: \ \  M_d(\P^n_a)\rightarrow N_d(\P^n_a)=W_d$$ is clearly a regular
map. The corresponding equivariant Euler class has the expression:

$$e_G(\P^n_a/W_d)=\prod_{j=0}^n\prod_{m=1}^{da_j}(a_jH-\lambda_j-m\alpha)$$with
$H$ the $T$-equivariant hyperplane class and $\lambda_j$'s the weights of the $T$-action.

Examples of singular toric varieties will be discussed again in our subsequent  paper, in which resolution of singularities will be used to reduce to the smooth case.
The above example was motivated by a question of Mazur, who 
suggested that the situation of counting rational curves in orbifolds is similar to certain Diophantine problem in number theory.

{\it Example 6:} For a general projective manifold $X$ embedded in
$\P(n)$, assume it is defined by a system of polynomial equations
$P(z^1, \cdots, z^n)=0$ where $z^j=(z^j_1, \cdots, z^j_{n_j})$ denotes the coordinate of $\P^{n_j}$. Assume the variety defined by the induced equation $P(f^1, \cdots,
f^n)=0$ in $N_d(\P(n))$ where $f^j=(f^j_1(w_0, w_1),\cdots,
f^j_{n_j}(w_0, w_1))$ is the tuple of polynomials, the coordinates for
the linear sigma model $N_d(\P^{n_j})$, is an orbifold.
Then we can take it to be our linear sigma model $W_d$. 
Note that $N_d(X)$ in Example 3 is embedded
inside this $W_d$. Very likely they are the same. 

Though we don't know whether this variety is an orbifold or not, 
it is clear that the fixed point compoenents in the above variety are given by $X_r$'s. 
In fact, we only need to assume that the localization formula holds on it. This is the case if the fixed point components embedded into $W_d$ as local complete intersection subvarieties. We conjecture that
 this is the case for any convex projective manifold.
Later, we will state a general conjectural Mirror Formula
in terms of this $W_d$.

\subsec{Regularity of the collapsing map}

For a toric manifold $X$, the following lemma show that $W_d$
described in Example 4 is 
a linear sigma model of $X$.
 
\lemma{ For toric manifold $X$, there is a regular extension 
$$\varphi: \ M_d(X)\rightarrow W_d$$ 
of the map $\varphi_o$ in Example 4 above.} 
\thmlab\Regularity

\let\po=\Po

\def\lra{\longrightarrow}
\def\cL{{\cal L}}
\def\cO{{\cal O}}
\def\cF{{\cal F}}
\def\cS{{\cal S}}

\def\cX{{\cal X}}
\def\sta{^{\ast}}
\def\mor{{\rm Mor}\,}

\def\cM{{\cal M}}

\proof  We simply follow the argument in \LLY, 
together with the construction in \Cox. We will define a morphism 
$\varphi:M_d(X)\ra W_d$.
Let $\cS$ be the category of all schemes of finite type (over $\bf C$)
and let
$$\cF: \cS\lra ({\rm Set})
$$
be the the contra-variant functor that send any $S\in \cS$ to the
set of families of stable morphisms
$$F: \cX\lra \po\times X\times S
$$
over $S$, where $\cX$ are families of connected arithmetic
genus 0 curves, modulo the obvious equivalence relation. Note that
$\cF$ is represented by the moduli stack $M_d(X) $.
Hence to define the morphism $\varphi$, 
it suffices to define a transformation
$$\Psi: \cF\lra \mor(-,W_d).
$$ 

We now define such a transformation. Let $S\in\cS$ and let
$\xi\in \cF(S)$ be represented by
$(\cX,F)$. We let $p_i$ be the
composite of $F$ with the $i$-th projection of
$\po\times X\times S$ and let $p_{ij}$ be the composite of $F$ with
the projection from $\po\times X\times S$ to the product of its
$i$-th and $j$-th components.
We consider the sheaf $p_2\sta\cO_{X}(L_j)$ on $\cX$ and
its direct image sheaf
$$\cL_{j,\xi}=p_{13\ast} p_2\sta\cO_{X}(L_j).
$$
Here the $L_j$ are the line bundles on $X$, as defined in Example 4.
As in \LLY, one can show that $\cL_{j,\xi}$ is flat in a standard way.


 For the same reasoning, the sheaves $\cL_{j,\xi}$
satisfy the following
base change property: let $\rho: T\to S$ be any base change and
let $\rho\sta(\xi)\in\cF(T)$ be the pull back of $\xi$. Then there is a
canonical isomorphism of sheaves of $\cO_T$-modules
\eqn\eqTwo{
 \cL_{j,\rho\sta(\xi)}\cong ({\bf 1}_{\po}\times \rho)\sta
\cL_{j,\xi}.
}

Since $\cL_{j,\xi}$ is flat over $S$, we can
define the determinant line bundle of $\cL_{j,\xi}$, denoted by $\det(\cL_{j,\xi})$ which is an invertible sheaf over $\po\times S$.
Using the Riemann-Roch theorem, one finds that its
degree along fibers over $S$ is $l_j=\bra c_1(L_j),d\ket$.
Furthermore, 
because $\cL_{j,\xi}$ has rank one, there is a canonical
homomorphism
\eqn\eqThree{
\cL_{j,\xi}\lra\det(\cL_{j,\xi}),
}
so that its kernel is the torsion subsheaf of $\cL_{j,\xi}$. 

Let $z_j$ be the $j$-th homogeneous
coordinate of $X$ (Example 4). 
Then $z_j$ is a section in $H^0(X,L_j)$. 
Its pull-back is a section of
$\cL_{j,\xi}$, which induces a section

$$\sigma_{j,\xi} \in H^0(S,\pi_{S\ast}\det(\cL_{j,\xi})).
$$
based on \eqThree.
Then after fixing an isomorphism
\eqn\eqOne{
\det(\cL_{j,\xi})\cong \pi_S\sta \cM\otimes\pi_{\po}\sta\cO_{\po}(l_j)
}
for some invertible sheaf $\cM$ of $\cO_S$-modules, 
where $l_j=\bra c_1(L_j), d\ket $.
We then obtain a section in 
$$ \pi_{S\ast}(\pi_{\po}\sta\cO_{\po}(l_j))\otimes_{\cO_S}\cM\equiv
H^0_{\po}(\cO_{\po}(l_j))\otimes_{{\bf C}}\cM.
$$ So $\sigma_{j, \xi}$ is determined up to certain constant $\lambda_j$ coming from $\cM$.

Now apply the above argument to each $L_j, \ j=1, \cdots, N$, we get $N$ sections $[\sigma_{1, \xi}, \cdots, \sigma_{N, \xi}]$. Let $w_0, w_1$ be the homogeneous coordinate of $\po$, we will write $\sigma_{j, \xi}=f_j(w_0, w_1)$ as a homogeneous polynomial of degree $l_j$. In this way we get a point in $\C^N(d)$. The constants $\lambda_j$ from ${\cal M}$ in choosing $\sigma_{j,\xi}$ must satisfy the relation $\prod_j \lambda_j^{\bra m,n_j\ket}=1$.
Here the $n_j$ are vectors in an integral lattice, which generate
the 1-cones in the defining fan of $X$, and
$m$ is any element in the dual lattice. (See \Cox.) 
 
For such $\lambda_j$'s we can then
 find an element $g$ in $T^m_C$ such that $[\lambda_1f_1, \cdots, \lambda_Nf_N]$ is transformed to $[f_1, \cdots, f_N]$ by $g$, therefore they represent the same point in $W_d$. In this way, after taking GIT quotient by $T_C^m$, the induced action on $[f_1(w_0, w_1),\cdots,f_N(w_0, w_1)]$ from the action on $\C^N$, for each map $(f, C)\in M_d(X)$, we have obtained canonically a point in $W_d$, therefore a morphism
$$\Psi(S): S\lra W_d
$$
that is independent of the isomorphisms \eqOne.
It follows from the base change property \eqTwo~ that the collection
$\Psi(S)$ defines a transformation
$${\bf \Psi}: \cF\lra \mor(-,W_d),
$$
thus defines the morphism $\varphi$ as desired.

The fact that $\varphi: M_d(X) \to W_d$ is $S^1\times T^N$-equivariant follows immediately from the fact that $\varphi$ is induced by
the transformation $\Psi$ of functors.
This completes the proof. $\Box$

For another proof of the above lemma we can proceed as follows. 
We use the notations as in the above Example 3. We show that the regularity of the collapsing map for $\P(n)$ induces the regularity of the collapsing map for $X$. For this we only need to prove that $N_d(X)$, the image  $\varphi(M_d(X))$ in $N_d(\P(n))$ of the collapsing map for $\P(n)$, lies in $W_d$. 

First, we show that $W_d$ lies in $N_d(\P(n))$. Note that both $W_d$ and $N_d(\P(n))$ are toric manifolds, and a Zariski open subset $W_d^o$ in $W_d$ is embedded $G$-equivariantly in $N_d(\P(n))$. 
Also the $G$-fixed points in $W_d$ are all in $X_r$, therefore in $\P(n)_r$ and $N_d(P(n))$. Any point in $W_d$ is in the closure of a generic $G_\C$ orbit in $W_d$ passing 
through two $G$ fixed points in $X_r$'s.
By the equivariance, this $G_\C$ orbit is also in $N_d(\P(n))$, therefore the closure of this orbit lies in $N_d(\P(n))$.

Second, we show that $N_d(X)$ lies in $W_d$. For this we note that $\varphi_o$ extends to $M_d(X)$, since it is actually the restriction of the 
corresponding map on $M_d(\P(n))$. Now by taking closure of the inclusion $\varphi(M^o_d(X))\subseteq W_d$, which is induced from the canonical identification, we get 

$$\varphi(M_d(X)) =\overline{\varphi(M^o_d(X))}\subseteq \overline{W_d}=W_d,$$
since $W_d$ is itself closed.

\newsec{The Gluing Identity}

Returning to the general case, 
we let $X$ be an admissible balloon manifold from now on.
In this section, we apply the functorial
 localization formula to the linear sigma
model. The argument used 
here is modelled on the one used in \LLY, 
except that the $T$-action is {\it not} used here.  
Thus all the results in this section hold for
manifolds without $T$ action.
We will have more to say about the mirror principle without $T$ action
later.

Recall that  we have
the commutative diagram:
$$\matrix{
F_r & {\br i_r\over \longrightarrow} & M_d(X)\cr
e\downarrow &  & \downarrow \varphi\cr
X_r & {\br j_r\over \longrightarrow} & W_d.}
$$
We also have the natural forgetting map
$\rho:M_{0,1}(d,X)\ra M_{0,0}(d,X)$,
and the projection map
$\pi: M_d(X)\ra M_{0,0}(d,X)$. Note 
that we have a commutative
diagram
$$\matrix{
M_d(X) &  &   \cr
\pi\downarrow & \nwarrow i_0 &  \cr
M_{0,0}(d,X) & {\br \rho\over\longleftarrow} & M_{0,1}(d,X).
}$$

Let $\varphi:M_d(X)\ra W_d,~e:F_r\ra X_r$ play the respective roles
of $f:X\ra Y,~g:F\ra E$ in the functorial localization formula.
Then it follows that

\lemma{Given any $G$-equivariant cohomology class $\omega$ on $M_d(X)$,
we have the following equality on $X_r$ for $0\preceq r\preceq d$:
$$\frac{j_r^*\varphi_!(\omega)}{e_G(X_r/W_d)}
=e_!\left(\frac{i_r^*(\omega)}{e_G(F_r/M_d(X))}\right).$$
}
\thmlab\SquareLemma

Actually this lemma may be viewed as an equivariant version of the so-called excess intersection formula of \Fulton, Theorem 6.3.

Let $L_r$ denote the line bundle on $ M_{0,
1}(r, X)$ whose fiber at $(f,C; x)$ is the tangent line at the marked point $x\in C$. Let $\pi_1$ denote the projection from $\P^1\times X$ to $\P^1$. 

The normal bundle of $F_r$ in $M_d(X)$ 
can be computed just as in \LLY.
For $r\neq 0, d$, we have 
$$N(F_r/M_d(X))=H^0(C_0, (\pi_1\circ f)^*T\P^1)+T_{x_1}C_0\otimes
L_r+T_{x_2}C_0\otimes L_{d-r}-A_{C_0}.$$ 
Here we have used the notations as
in \LLY: a point $(f_1, C_1, x_1)$ in $ M_{0,
1}(r, X)$ and a point $(f_2, C_2, x_2)$ in $M_{0,
1}(d-r, X)$ is glued to $C_0\simeq \P^1$ at $0$ and $\infty$ respectively to get the point $(f, C)$ in $M_d(X) $ with $C\simeq C_1 \cup C_0\cup C_2$. Since $x_1$ and $x_2$ are mapped to the same point in $X$ under the  projection $\pi_2: \ \P^1\times X\rightarrow X$, so this point can be considered as a point in $F_r$ by gluing together $(f_1, C_1, x_1)$ and $(f_2, C_2, x_2)$ at the marked points.
Similarly, for $r=0, d$, we
have 
$$N(F_0/M_d(X))=H^0(C_0, (\pi_1\circ f)^*T\P^1)+T_{x_1}C_0\otimes L_d-A_{C_0}$$
and 
$$N(F_d/M_d(X))=H^0(C_0,(\pi_1\circ  f)^*T\P^1)+T_{x_2}C_0\otimes L_{d}-A_{C_0}.$$
In the above $H^0(C_0, (\pi_1\circ f)^*T\P^1)$ corresponds to the deformation of $C_0$;
 $T_{x_1}C_0\otimes L_r$ and $T_{x_2}C_0\otimes L_{d-r}$ 
correspond respectively to the deformations of the nodal points $x_1$ and $x_2$;  $A_{C_0}$ denotes the automorphism group  to be quotiented out.

The equivariant Euler classes of the normal bundles above
are computed  as in \LLY, to which we refer the readers for details.
For $r\neq 0, d$,
the equivariant Euler classes are:
$$e_G(F_r/M_d(X))=
-\alpha(-\alpha+c_1(L_{d-r})) \cdot
\alpha(\alpha+c_1(L_r))
$$ 
where the two factors on the right hand side
are pullbacked to $F_r$ from $M_{0,1}(d-r,X)$, $M_{0,1}(r,X)$
respectively.
For $r=0,d$, we have
$$e_G(F_0/M_d(X))=-\alpha(-\alpha+c_1(L_d)),\
e_G(F_d/M_d(X))=\alpha(\alpha+c_1(L_d))$$
respectively. Combining this with the preceding lemma, we get
the following equality on $X=X_0$:
$$\frac{j_0^*\varphi_!(\omega)}{e_G(X_0/W_d)}
=ev_!\left(\frac{i_0^*(\omega)}{\alpha(\alpha-c_1(L_d))}\right).$$
Here we have dropped the subscript from $ev_d$.
In particular, if $\psi$ is a class on $M_{0,0}(d,X)$,
then for $\omega=\pi^*\psi$, we get
$i_0^*(\omega)=i_0^*(\pi^*\psi)=\rho^*\psi$.
This yields

\lemma{Given any $T$-equivariant cohomology class 
$\psi$ on $M_{0,0}(d,X)$,
we have the following equality on $X$:
$$\frac{j_0^*\varphi_!(\pi^*\psi)}{e_G(X_0/W_d)}
=ev_!\left(\frac{\rho^*\psi}{\alpha(\alpha-c_1(L_d))}\right).$$
}
\thmlab\AlphaSquare

\lemma{For $0\preceq r\preceq d$, we have
the following equality on $X$:
$$e_G(X_r/W_d)=\overline{e_G(X_0/W_r)} e_G(X_0/W_{d-r}).$$
In particular, we have 
$$e_G(X_d/W_d)=\overline{e_G(X_0/W_d)}.$$
}
\thmlab\eee
\proof
Consider the commutative diagram
\eqn\FiberSquare{\matrix{
F_r & {\br \Delta_0\over\longrightarrow}& 
M_{0,1}(r,X)\times M_{0,1}(d-r,X)\cr
e\downarrow &  &~~~~~~~\downarrow ev_r\times ev_{d-r}\cr
X & {\br \Delta\over\longrightarrow} & X\times X
}}
where $\Delta$ is the diagonal map,
and $\Delta_0$ is the inclusion
induced by $\Delta$. In particular, by definition we
 have $(ev_r\times ev_{d-r})^*\Delta_!(1)=(\Delta_0)_!(1)$. 
So we have
\eqn\dumb{\Delta^*(ev_r\times ev_{d-r})_!(\omega)
=e_!\Delta_0^*(\omega)
}
for any class $\omega$ on
$M_{0,1}(r,X)\times M_{0,1}(d-r,X)$.
Now put $\omega={1\over\alpha(\alpha+c_1(L_r))}\times
{1\over\alpha(\alpha -c_1(L_{d-r}))}$.
Then \dumb~
becomes
$$
(ev_r)_!{1\over \overline{e_G(F_0/M_r(X))} }
\cdot (ev_{d-r})_!{1\over e_G(F_0/M_{d-r}(X)) }
=e_!{1\over e_G(F_r/M_d(X))}.
$$ 

Since $\varphi:M_d(X)\ra W_d$ is an isomorphism on a 
Zariski open set, we see that $\varphi_!(1)=1$. In fact, by 
Prop. (5.3.3) \AP,  $\varphi_!$ preserves degree because
since $M_d(X)$ and $W_d$ have the same dimension.
So $\varphi_!(1)\in H^0_G(W_d)$ is a constant. By 
restricting to the Zariski open set on which $\varphi$ is an isomorphism, we 
find $\varphi_!(1)=1$.

By taking $\psi=1$ in the preceding lemma, we get
$$(ev_r)_!{1\over e_G(F_0/M_r(X))}={1\over e_G(X_0/W_r)},~~~
(ev_{d-r})_!{1\over e_G(F_0/M_{d-r}(X))}={1\over e_G(X_0/W_{d-r})}.
$$
By taking $\omega=1$ in Lemma \SquareLemma, we get
$$e_!{1\over e_G(
F_r/M_d(X))}={1\over e_G(X_r/W_d)}.$$
Combining the last four equations yields our assertion.
$\Box$

Fix a $T$-equivariant multiplicative class $b_T$.
Fix a $T$-equivariant bundle of the form $V=V^+\oplus V^-$,
where $V^\pm$ are respectively the convex/concave bundles.
(cf. \LLY.)
We call such a $V$ a mixed bundle.
We assume that
$$\Omega:={b_T(V^+)\over b_T(V^-)}$$
is a well-defined invertible class on $X$.
By convention, if $V=V^\pm$ is purely convex/concave,
then $\Omega=b_T(V^\pm)^{\pm1}$.
Recall that the bundle $V\ra X$ induces the bundles
$$V_d\ra M_{0,0}(d,X),~~
U_d\ra M_{0,1}(d,X),~~
\cU_d\ra M_d(X).$$
Moreover, they are related by $U_d=\rho^* V_d$, $\cU_d=\pi^* V_d$,
Throughout this section, we denote
$$
Q:~~Q_d:=\varphi_!(\pi^* b_T(V_d)).
$$
If $\omega$ is a class on $W_d$,
we write
$$
i_r^*\omega^v:={j_r^*\omega\over
e_G(X_r/W_d)}
$$
which is a class on $X=X_r$.

\lemma{For $0\preceq r \preceq  d$, 
$$\Omega~i_r^*Q_d^v=\overline{i_0^*Q_r^v}~
i_0^*Q_{d-r}^v.$$
}
\proof
For simplicity,
 let's consider the case $V=V^+$.
The general case is entirely analogous.

Recall that a point $(f,C)$ in $F_r\subset M_d$ comes
from gluing together
a pair of stable maps $(f_1,C_1,x_1),(f_2,C_2,x_2)$
with $f_1(x_1)=f_2(x_2)=p \in X$. From this, we get an exact sequence over $C$:
$$0\rightarrow f^*V\rightarrow f_1^*V\oplus f_2^*V\rightarrow
V|_p\rightarrow 0.$$ 
Passing to cohomology, we have   
$$0\rightarrow H^0(C, f^*V)\rightarrow H^0(C_1, f_1^*V)\oplus  H^0(C_2,
f_2^*V)
\rightarrow V|_p\rightarrow 0.$$
Hence we obtain an exact sequence of bundles
on $F_r$:
$$0\rightarrow i_r^*\cU_d\rightarrow U_r'\oplus U_{d-r}'\rightarrow
e^*V\rightarrow 0.$$
Here $i_r^*\cU_d$ is the restriction to $F_r$ of the bundle
 $\cU_d\ra M_d(X)$. And $U_r'$ is the pullback of the bundle
$U_r\ra M_{0,1}(d,X)$, and similarly for $U_{d-r}'$.
Taking the multiplicative characteristic class $b_T$, we get
the identity on $F_r$:
$$e^*b_T(V)b_T(i^*_r\cU_d)=b_T(U_r')b_T(U_{d-r}').$$ 
This is what we call
the {\it gluing identity}.

Now put
$$\omega={b_T(U_r)\over e_G(F_r/M_r(X))}\times
{b_T(U_{d-r})\over e_G(F_0/M_{d-r}(X))}.$$
From the commutative diagram \FiberSquare, we have the
identity:
$$
\Delta^*(ev_r\times ev_{d-r})_!(\omega)
=e_!\Delta_0^*(\omega).
$$
On one hand is
$$\eqalign{
\Delta^*(ev_r\times ev_{d-r})_!(\omega)
&=(ev_r)_!{b_T(U_r)\over e_G(F_r/M_r(X))}~\cdot~
(ev_{d-r})_!{b_T(U_{d-r})\over e_G(F_0/M_{d-r}(X))}\cr
&=(ev_r)_!{\rho^*b_T(V_r)\over e_G(F_r/M_r(X))}~\cdot~
(ev_{d-r})_!{\rho^*b_T(V_{d-r})\over e_G(F_0/M_{d-r}(X))}\cr
&=\overline{i_0^*Q_r^v}~i_0^*Q_{d-r}^v,
}
$$
the last equality being a consequence of Lemma \AlphaSquare.
On the other hand, applying the gluing identity, we have 
$$\eqalign{
e_!\Delta_0^*(\omega)
&=e_!\left({b_T(U_r')\over\alpha(\alpha+c_1(L_r))}~
{b_T(U_{d-r}')\over\alpha(\alpha-c_1(L_{d-r}))}\right)\cr
&=e_!{e^*b_T(V)i^*_rb_T(\cU_d)\over e_G(F_r/M_d(X))}\cr
&=b_T(V)~e_!{i^*_rb_T(\cU_d)\over e_G(F_r/M_d(X))}\cr
&=b_T(V) i_r^*Q_d^v,
}
$$
the last equality being a consequence of Lemma \SquareLemma.
This proves our assertion.  $\Box$

\remark{\item{(a)} If we take $V$ to be the trivial line bundle,
and $b_T$ to be the total Chern class, then the preceding
lemma reduces to Lemma \eee. 
\item{(b)} All the lemmas in this section, in fact, holds
for a general projective manifold $X$ without $T$-action,
provided that we still have the $S^1$-equivariant map
$\varphi:M_d(X)\ra W_d$, with properties 1.-2. stated
in section \SigmaModels. All $G$-equivariant classes
above are then replaced by their $S^1$-equivariant counterparts.
}

\newsec{Euler Data}

{\it Notations:} We denote by $\kappa_i$ the $G$-equivariant
class on $W_d$ with the property that $j_r^*\kappa_i=H_i+\bra
H_i,r\ket\alpha$. 
By the localization theorem, 
$\kappa_i$ is determined by these restriction conditions,
and is a class in the localized equivariant cohomology of $W_d$.
More generally a class $\phi\in H^2_T(X)$ has
a $G$-equivariant lift $\hat\phi\in H_G^2(W_d)$
determined by $j_r^*\hat\phi=\phi+\bra\phi,r\ket\alpha$.
We denote by $\bra H^2_T(X)\ket$
the ring generated by $H^2_T(X)$, and
by $R_d$ the ring generated by their lifts $\hat\phi$.
We put $\cR=\Q(\cT^*)[\alpha]$, where $\Q(\cT^*)$ is the
rational function field on the Lie algebra of $T$.
For convenience, we introduce the notations
$$\eqalign{
\kappa\cdot\zeta=\kappa_\zeta&:=\kappa_1\zeta_1+\cdots+
\kappa_m\zeta_m\cr
i_r^*\omega^v&:={j_r^*\omega\over
e_G(X_r/W_d)}
}$$
where $\omega$ is a class on $W_d$.

It is often necessary to work over a larger
field than
$\C$ for coefficients of cohomology groups. 
For example when we consider 
the case of the equivariant Chern polynomial $c_T$,
a formal variable $x$ is introduced. In this
case we replace everywhere the scalars $\C$ by $\C(x)$.
This will be implicit in all of the discussion below.

Recall the localization formula: 
$$\int_{W_d}\omega
=\sum_{0\preceq r\preceq d}\int_{X}
{j^*_r(\omega)\over e_G(X_r/W_d)}.$$
We shall often apply the following version:
$$\int_{W_d}\omega~e^{\kappa_\zeta}
=\sum_{0\preceq r\preceq d}\int_{X}
i^*_r\omega^v~ e^{H_\zeta+\bra H_\zeta,r\ket\alpha}.$$

\definition{Fix an invertible class
$\Omega\in H_T^*(X)^{-1}$.
A list $P:~P_d\in H_G^*(W_d)^{-1}$, $d\succ0$, is
a  $\Omega$-Euler data if on $X$,
$$\Omega~ i_r^*P_d^v=
\overline{i_0^*P_r^v}~i^*_0P_{d-r}^v$$
(called Euler data identity)
for all $r\preceq d$, and
the $\int_{W_d}P_d\cdot\omega$
are polynomial in $\alpha$ for all $\omega\in R_d$.
By convention we set $P_0=\Omega$.
}

{\it Example 0.} In the last section
we have  proved,
using the {\it gluing identity}, 
that the data $Q:~~Q_d=\varphi_!(\pi^*b_T(V_d))$
associated with a mixed bundle $V$
and a multiplicative class $b_T$ satisfies
the Euler data identity. This indicates that
the gluing identity is really the
geometric origin of Euler data. This is what motivates
our definition of Euler data. Note that since $Q_d$ is the
equivariant push-forward of a class in $H_G^*(M_d(X))$,
the polynomial condition on $Q$ is automatic. This condition 
will be needed 
when mirror transformation is discussed.

{\it Example 1.} Let $L$ be any equivariant line bundle with 
$c_1(L)\geq0$. Let $\hat L$ be the $G$-equivariant lift of
$c_1(L)$. 
$$P_d=\prod_{k=0}^{\bra c_1(L),d\ket}(\hat L-k\alpha)$$
is an $\Omega$-Euler data where $\Omega=c_1(L)$.

{\it Example 2.} Let $L$ be any equivariant line bundle with 
$c_1(L)< 0$. Let $\hat L$ be the $G$-equivariant lift of
$c_1(L)$. 
$$P_d=\prod_{k=1}^{-\bra c_1(L),d\ket-1}(\hat L+k\alpha)$$
is an $\Omega$-Euler data where  $\Omega=c_1(L)^{-1}$. 

{\it Example 3.} If $P,P'$ are $\Omega,
\Omega'$-Euler data
respectively, then $P\cdot P'$ is a 
$\Omega\Omega'$-Euler data
as shown in \LLY.

{\it Example 4.} Let $L$ be as in Example 1, and $x$
be a formal variable. Then
$$P_d=\prod_{k=0}^{\bra c_1(L),d\ket}(x+\hat L-k\alpha)$$
is an $\Omega$-Euler data where $\Omega=c_T(L)$ denotes the
Chern polynomial.

\bs
In each of Examples 1-4 above, the Euler data 
identity  follows immediately from the algebraic
identity $\Omega~j_r^*P_d=\overline{j_0^* P_r}~
j_0^*P_{d-r}$, and
Lemma \eee.

Strictly speaking, in the examples above, we must
require that $c_1(L)$ be an invertible class. This
requirement can be easily met by twisting $L$ by
a trivial line bundle on which $T$ acts by a suitable
weight. In the end, we will only be interested in
the nonequivariant limit of an Euler data. Thus
the choice of twisting is of no consequence at
the end. Alternatively, we can consider the
Chern polynomial or
the total Chern class (which is automatically invertible)
instead of the first Chern class. 

\subsec{An algebraic property}

Let $\cS$ denotes
the set of sequences
$B:~~B_d\in H_G^*(X)^{-1},~~d\succeq0$.
By convention, we set $B_0=\Omega$.

\definition{
Given any $B\in\cS$, define the formal series
$$
HG[B](t):=
e^{-H\cdot t/\alpha}\left( \Omega+ 
\sum_{d\succ 0}
B_d~e^{d\cdot t}\right).
$$
}
Note that $e^{H\cdot t/\alpha}HG[B](t)$ takes value in the ring
$H_G(X)^{-1}[[K^\vee]]$. (Notations: if $R$ is a ring, 
hen
$\cR[[K^\vee]]:=\{\sum_{d\in\Lambda} a_d e^{d\cdot t}|a_d\in\cR\}$.
We use the notations $e^{d\cdot t}=e^{\bra H_t,d\ket}$
interchangeably.)

Let $P$ be an Euler data, and let $B$ be the
list with $B_d:=i_0^*P_d^v$.
By the localization formula and the Euler data identity,
we have
$$\eqalign{
\int_{W_d} P_d~ e^{\kappa_\zeta}
&=\sum_{r\preceq d}\int_X
i^*_r P_d^v~ e^{H_\zeta+\bra H_\zeta,r\ket\alpha}\cr
&=\sum_{r\preceq d}e^{-d\cdot\tau}\int_X\Omega^{-1}
\overline{\left[e^{-H_t/\alpha}i_0^* P_r^v~e^{r\cdot t}\right]}~
\left[e^{-H_\tau/\alpha}i_0^*P_{d-r}^v~e^{(d-r)\cdot \tau}\right].
}
$$
Here $t=\zeta\alpha+\tau$.
Note that $\bar\zeta=-\zeta$, $\bar\alpha=-\alpha$,
and all other variables are invariant under
the ``bar'' operation.
Now multiply both sides by $e^{d\cdot\tau}$ and sum over
$d\in K^\vee$, we get the formula:
\eqn\THREE{
\sum_d e^{d\cdot\tau}\int_{W_d} P_d~ e^{\kappa_\zeta}
=\int_X\Omega^{-1}~
\overline{HG[B](\zeta\alpha+\tau)}~HG[B](\tau).
}
By definition, the coefficient of
$e^{d\cdot\tau}$ on the right hand side is a power series
in $\zeta$ with coefficients which are polynomial in
$\alpha$, ie. the series lies in $\cR[[e^\tau,\zeta]]$.

Conversely, given $B\in\cS$ such that 
$$\int_X\Omega^{-1}~
\overline{HG[B](\zeta\alpha+\tau)}
~HG[B](\tau)\in\cR[[e^\tau,\zeta]],$$
there exists a unique Euler data $P:~P_d$
satisfying \THREE. Namely, $P_d$ is defined by
the conditions
$$j_r^*P_d=\Omega^{-1}e_G(X_r/W_d)~
\overline{B_r}~B_{d-r}.$$

Thus an Euler data $P$ gives  rise to
a list $B\in\cS$ in a canonical way. 
Abusing the terminology, {\it we shall call
such a $B$ an Euler data}.

\newsec{Linking and Uniqueness}

\lemma{Let $\omega\in H^*_T(X)^{-1}(\alpha)$. Suppose that
\item{(a)} for each $q\in X^T$,
$\omega_q(\alpha):=\omega(\alpha)|_q$ is a Laurent polynomial
in $\alpha$ with $deg_\alpha\omega(\alpha)\leq -2$;
\item{(b)} the power series in $\zeta$:
$\int_X\left(\omega(\alpha) e^{H_\zeta} + \omega(-\alpha) 
e^{H_\zeta+\bra H_\zeta,d\ket~\alpha}\right)$
has coefficients
which are polynomial in $\alpha$.
\item{}
Then $\omega=0$.
}
\proof
Suppose $\omega\neq0$, and we will get a contradiction.
By assumption (a), we can write
$${\omega_q(\alpha)\over e_T(q/X)}
=a_q\alpha^{-k}+b_q\alpha^{-k+1}+\cdots$$
which is a {\it finite} sum, with $a_q$
independent of $\alpha$ and  $k\geq2$. By supposition,
not all the $a_q$ are zero.
By localization, we get
$$\eqalign{
&\int_X\left(\omega(\alpha) e^{H_\zeta} + \omega(-\alpha) 
e^{H_\zeta+\bra H_\zeta,d\ket~\alpha}\right)\cr
&=\sum_q\left((a_q\alpha^{-k}+b_q\alpha^{-k+1}+\cdots)
e^{H_\zeta(q)}
+(a_q(-\alpha)^{-k}+b_q(-\alpha)^{-k+1}+\cdots)
e^{H_\zeta(q)+\bra H_\zeta,d\ket\alpha}\right).}$$
By assumption (b), order by order in $\zeta$,
this expression is polynomial in $\alpha$. Since
$k\geq2$, the polar term with $\alpha^{-k}$ must vanish,
and so 
$$\sum_q a_q e^{H_\zeta(q)}(1+(-1)^k)=0.$$
Since not all $a_q$ are zero and the functions $e^{H_\zeta(q)}$
are linearly independent,
it follows that $k$ is odd.
Now the coefficient of $\alpha^{-k+1}$ becomes
$$\sum_q e^{H_\zeta(q)}
(2b_q-a_q  \bra H_\zeta,d\ket)=0.$$
Again by linear independence of the exponential functions,
it follows that $a_q=0=b_q$ for all $q$, which
is a contradiction.  $\Box$

\lemma{Suppose $A,B$ are Euler data with
$A_r=B_r$ for all $r\prec d$. 
Suppose that the
$(A_d-B_d)|_q$, $q\in X^T$, 
are Laurent polynomial in $\alpha$. 
Suppose also that $deg_\alpha (A_d-B_d)
\leq -2$.
Then $A_d=B_d$.}
\thmlab\LaurentUniqueness
\proof
It suffices to show that
$$\omega(\alpha):=A_d-B_d$$
has property (b) of the preceding lemma.

Let $A_d'$ be the coefficient of $e^{d\cdot\tau}$
in the series 
$$\int_X\Omega^{-1}
\overline{HG[A](\zeta\alpha+\tau)}
HG[A](\tau).$$
Likewise for the $B_d'$. Since $A,B$ are Euler data,
$A_d',B_d'$ are power series in $\zeta$
with coefficients which are polynomial in $\alpha$.
Explicitly,
$$
A_d'=\sum_{r\preceq d}\int_X\Omega^{-1}
e^{H_\zeta}A_r e^{r\cdot \zeta\alpha}~
A_{d-r}
$$
and likewise for the $B_d'$. Using that $A_r=B_r$, $r\prec d$,
and that $A_0=B_0=\Omega$, we see that
$A_d'-B_d'$ is a sum over $r$ with only
two surviving terms, corresponding to $r=0,d$.
That is,
$$
A_d'-B_d'=
\int_X\left(\omega(\alpha) e^{H_\zeta} + \omega(-\alpha) 
e^{H_\zeta+\bra H_\zeta,d\ket~\alpha}\right).
$$
Since both $A_d',B_d'$ have coefficients
which are polynomial in $\alpha$,
this shows that the class $\omega(\alpha)$
has property (b) of the preceding lemma.
$\Box$

\definition{Two Euler data $A,B$ are linked if
for every balloon $pq$ in $X$ and
every $d=\delta[pq]\succ0$,
$$(A_d-B_d)|_q$$
is regular
at $\alpha={\lambda\over\delta}$
where $\lambda$ is the weight on the tangent line $T_q(pq)$.
}

Suppose $A,B$ both come from Euler data $Q,P$ respectively,
ie. $A_d=i_0^*Q_d^v$ and $B_d=i_0^*P_d^v$. Suppose also
that
\eqn\EulerLinking{
j_0^*(P_d)|_q=j_0^*(Q_d)|_q~~~~at~\alpha=\lambda/\delta.
}
whenever $d=\delta[pq]\succ0$ as above. 
Recall that 
$\alpha=\lambda/\delta$ is at worst a simple pole
of $1/e_G(X_0/W_d)|_q$. It follows that
$(A_d-B_d)|_q$ is regular at this value.
This shows that the conditions \EulerLinking~
guarantee that $A,B$ are linked.

\theorem{Suppose $A,B$ are linked Euler data satisfying
the following properties: for $d\succ0$,
\item{(i)} If $q\in X^T$,
the only possible poles of
$(A_d-B_d)|_q$
are scalar multiples of a weight on $T_qX$.
\item{(ii)} 
$deg_\alpha (A_d-B_d)\leq -2$.
\item{} Then $A=B$.}
\thmlab\UniquenessTheorem
\proof
We will prove, by induction, the assertion 
that $A_d=B_d$ for all $d$. If $d=0$, there is
nothing to prove. Suppose the assertion holds for
all $r\prec d$. Set
$\omega_q(\alpha):=(A_d-B_d)|_q$ as before.
We will show, under assumption (i),
that the $\omega_q(\alpha)$
are Laurent polynomial in $\alpha$.
It follows then, from the preceding lemma and assumption (ii), that
$A_d=B_d$.

Let $\lambda\in\cT^*-0$.
We will show that each $\omega_q(\alpha)$ is regular
at $\alpha=\lambda$.
Recall the power series in $\zeta$:  $A_d',B_d'$,
with coefficients polynomial in $\alpha$
as in the preceding proof.
Thus for any integers $k,l\geq0$,
$$Res_{\alpha=\lambda}\left((\alpha-\lambda)^k(\alpha+\lambda)^l
(A_d'-B_d')\right)=0.$$
Also recall that
$$\eqalign{
A_d'-B_d'
&=\int_X\left(\omega(\alpha) e^{H_\zeta} + \omega(-\alpha) 
e^{H_\zeta+\bra H_\zeta,d\ket~\alpha}\right)\cr
&=\sum_{q\in X^T}{1\over e_T(q/X)}\left(\omega_q(\alpha)
 ~e^{H_\zeta(q)}
+\omega_q(-\alpha)
~e^{H_\zeta(q)+\bra H_\zeta,d\ket\alpha}
\right).
}
$$
From the preceding two equations, we get
\eqn\dumb{\eqalign{
0&=\sum_{q\in X^T}{1\over e_T(q/X)}\left(
e^{H_\zeta(q)}~
Res_{\alpha=\lambda}(\alpha-\lambda)^k(\alpha+\lambda)^l\omega_q(\alpha)
\right.\cr
&\left.~~~~~~~+e^{H_\zeta(q)+\bra H_\zeta,d\ket \lambda}
Res_{\alpha=\lambda}(\alpha-\lambda)^k(\alpha+\lambda)^l\omega_q(-\alpha)
\right).
}
}
If 
$Res_{\alpha=\lambda}(\alpha-\lambda)^k(\alpha+\lambda)^l\omega_q(-\alpha)=0$
for all $q$, then the preceding equation shows that
$Res_{\alpha=\lambda}(\alpha-\lambda)^k(\alpha+\lambda)^l\omega_q(\alpha)=0$
for all $q$, because 
the vectors $H_\zeta(q)$ are
distinct.
Similarly if
$Res_{\alpha=\lambda}(\alpha-\lambda)^k(\alpha+\lambda)^l\omega_q(\alpha)=0$
for all $q$, then  we have
$Res_{\alpha=\lambda}(\alpha-\lambda)^k(\alpha+\lambda)^l\omega_q(-\alpha)=0$
for all $q$. In either case, we conclude that each
$\omega_q(\alpha)$ is regular at $\alpha=\lambda$.
So if $\alpha=\lambda$ is a pole of a $\omega_q(\alpha)$, then
we necessarily have
$Res_{\alpha=\lambda}(\alpha-\lambda)^k(\alpha+\lambda)^l\omega_q(\alpha)\neq0$ 
{\it and}
$Res_{\alpha=\lambda}(\alpha-\lambda)^k(\alpha+\lambda)^l\omega_p(-\alpha)\neq0$ 
for some  $p,q$ and some $l,k$, such that 
$$H_\zeta(q)=H_\zeta(p)+\bra H_\zeta,d\ket\lambda,$$
to ensure cancellation of the exponential functions in \dumb.
Note that since $d\succ0$ and $\lambda\neq0$, we have $p\neq q$.
By our assumption (i),
the pole $\alpha=\lambda$ of $\omega_q(\alpha)$
must be of the form $\lambda={\lambda'\over \delta}\neq0$ 
for some weight $\lambda'$ on $T_qX$, and some scalar $\delta\neq0$. 
By Lemma \BalloonLemma,  $p,q$ must be joined by
a balloon, 
$d=\delta[pq]$,
and $\lambda'$ is the weight on the tangent line $T_q(pq)$.
Thus if $d$ is not a multiple of $[pq]$, then
we have shown that $\omega_q(\alpha)$ is regular away from
$\alpha=0$.

Now suppose that $d=\delta
[pq]$, and consider the only possible pole of $\omega_q(\alpha)$
at $\alpha={\lambda'\over\delta}\neq0$, as above.
By hypothesis, $A,B$ are linked.
But this means that $\omega_q(\alpha)$ is regular
at $\alpha={\lambda'\over\delta}\neq0$.
$\Box$

\remark{In our applications later, the situation is
better then the conditions (i)-(ii) demand.
We will have two Euler data $A,B$ such that
$A_d,B_d$ separately, rather than just $A_d-B_d$,
will satisfy both conditions (i)-(ii) at the
outset. In this situation, to prove that $A=B$,
it suffices to prove that they are linked.}
\thmlab\UniquenessRemark

\newsec{Mirror Transformations}

Throughout this section, we fix an invertible class $\Omega$
on $X$, and
will denote by $\cA$ the set of
$\Omega$-Euler data.

\definition{A map $\mu:\cA\ra\cA$ 
is called a mirror transformation if it preserves
linking. In other words, $\mu(A)$ and $A$ are linked
for any $A\in\cA$.
We call $\mu(A)$ a mirror transform of $A$.}

We now consider a construction of mirror transformations,
as motivated by the classic example of \CDGP. Consider a
transformation $\mu:\cS\ra \cS$, 
$B\ra\tilde B$, of the type
\eqn\GeneralTransformation{
\tilde B_d
=B_d+\sum_{r\prec d} a_{d,r} B_r
}
where the $a_{d,r}\in H_G^*(X)^{-1}$ are a given
set of coefficients.
This transformation is obviously invertible, and preserves
$B_0=\Omega$. 

\lemma{Suppose that $B,\tilde B$ are both
Euler data. 
Let $d=\delta[pq]\succ0$ for
some balloon $pq$ in $X$.
Suppose that the coefficients in 
\GeneralTransformation~ are such that
their restrictions $a_{d,r}(q)$, 
$r\prec d$, to the fixed point $q$ are regular
at $\alpha=\lambda/\delta$ where $\lambda$ is
the weight on $T_q(pq)$. Then
$(\tilde B_d-B_d)|_q$ is regular at 
$\alpha=\lambda/\delta$.
}
\thmlab\RegularityLemma
\proof
From \GeneralTransformation,
it suffices to show that  the functions
$B_r|_q$, $0\prec r\prec d$, are regular
at $\alpha=\lambda/\delta$. Suppose the contrary that
some $B_r|_q =0$ has a pole of order $k+1$ there.
Since $B$ is a Euler data, we know that
$$
B_r':=\sum_{s\preceq r}\int_X\Omega^{-1}
e^{H_\zeta}B_s e^{\bra H_\zeta,s\ket\alpha}~
B_{r-s}
$$
is a power series in $\zeta$ with coefficients
polynomial in $\alpha$.

By the localization formula,
$$
B_r'=\sum_{s\preceq r}\sum_{o\in X^T}
{1\over e_T(o/X)} 
e^{H_\zeta(o)+\bra H_\zeta,s\ket\alpha}
\Omega(o)^{-1}B_s(o)B_{r-s}(o).
$$
Now multiply both sides by $(\alpha-\lambda/\delta)^k$ and
take residue at $\alpha=\lambda/\delta$. We get 
$$0=\sum_{s\preceq r}\sum_{o\in X^T}
{1\over e_T(o/X)} 
e^{H_\zeta(o)+\bra H_\zeta,s\ket\lambda/\delta}
Res_{\alpha=\lambda/\delta}(\alpha-\lambda/\delta)^k
\Omega(o)^{-1}B_s(o)B_{r-s}(o).
$$
By assumption, the summand above with $s=0, o=q$ is
nonzero. Observe that this term has an exponential factor
$e^{H_\zeta(q)}$. Thus in order to cancel this term,
any other term contributing to this cancellation must
have an identical exponential factor. This means that
$$H_\zeta(q)=H_\zeta(o)+\bra H_\zeta,s\ket\lambda/\delta$$ 
for some $s$ with $s\preceq r$, and some $o\in X^T$.
By Lemma \BalloonLemma, this implies that $s=\delta [pq]$,
contradicting that $s\preceq r\prec d$. 
$\Box$

\definition{The transformation \GeneralTransformation~
is said to have the regularity property if
for every balloon $pq$ in $X$ and $d=\delta[pq]$, 
the coefficients
are such that their restrictions $a_{d,r}(q)$, $r\prec d$,
are regular 
at $\alpha=\lambda/\delta$ where $\lambda$ is
the weight on $T_q(pq)$.
} 
Thus the preceding lemma says that transformation
\GeneralTransformation~ having the regularity property
preserves linking.

Again, motivated by \CDGP~ and \HLYII, we consider the following
special types of transformations.
Given a power series $f\in\cR[[K^\vee]]$ with no constant term,
we have an invertible transformation $\mu_f:\cS\ra\cS$,
$B\mapsto\tilde B$, such that
$$e^{f/\alpha}~HG[B](t)
=HG[\tilde B](t).$$
In fact, we have
$$\tilde B_d=B_d
+\sum_{r\prec d} f_{d-r} B_r$$
where 
$e^{f/\alpha}=\sum_{s\succeq0}f_s e^{s\cdot t}$, 
$f_s\in\cR[\alpha^{-1}]$.
This is clearly a transformation of type \GeneralTransformation
~ having the regularity property. (In fact, 
all the coefficients $f_{d-r}$ 
are regular away from $\alpha=0$.)

Given power series $g=(g_1,..,g_m)$,
$g_j\in\cR[[K^\vee]]$ with no constant term,
we have an invertible transformation $\nu_g:\cS\ra\cS$,
$B\mapsto\tilde B$, such that
$$HG[B](t+g)
=HG[\tilde B](t).$$
In fact since
$$
HG[B](t+g)
=e^{-H\cdot t/\alpha} e^{-H\cdot g/\alpha}\sum_{d\succeq0} 
B_d~ e^{d\cdot t} e^{d\cdot g},
$$
if we write $e^{d\cdot g}=\sum_{s\succeq0} g_{d,s}e^{s\cdot t}$,
$g_{d,s}\in\cR$ and
$e^{-H\cdot g/\alpha}=\sum_{s\succeq0} \hat g_s e^{s\cdot t}$,
$\hat g_s\in\cR[H/\alpha]$, then 
$$\tilde B_d
=B_d+\sum_{r\prec d} a_{d,r} B_r$$
where the $a_{d,r}\in H_G^*(X)^{-1}$ are quadratic
expressions in the $g,\hat g$.
Thus we obtain another transformation $\cS\ra\cS$
of type \GeneralTransformation, again having the regularity property.

\theorem{The transformations $\mu_f,\nu_g:B\mapsto \tilde B$ above
each defines a mirror transformation. That is,
if $B$ is a Euler data
then $\mu_f(B)$ and $\nu_g(B)$ are both
Euler data linked to $B$.
}
\thmlab\MirrorTransform
\proof 
Let $B$ be a given Euler data.
We have seen that the preceding lemma
guarantees that $\mu_f(B)$, $\nu_g(B)$
are linked to $B$.
So it suffices to show that they are
Euler data.

First case: set $\tilde B=\nu_g(B)$, ie.
\eqn\dumb{
HG[\tilde B](t)=HG[B](t+g(e^t)).
}
(Here $e^t$ means the variables $(e^{t_1},..,e^{t_m})$.)
Set $t=\zeta\alpha+\tau$, $q_i=e^{\tau_i}$.
On the one hand, we have
\eqn\THREE{
\int_X\Omega^{-1}~
\overline{HG[B](t)}~ HG[B](\tau)
=\sum_{d,m} 
q^d~\zeta^m~B'_{d,m}(\alpha)
}
for some $B'_{d,m}\in\cR$. 
Now compare
$$\eqalign{
(*)~~~\int_X\Omega^{-1}~
\overline{HG[B](t)}~ HG[B](\tau)
&=\int_X\Omega^{-1}e^{H_\zeta} 
\sum \bar B_d e^{d\cdot\tau} e^{d\cdot\zeta\alpha}\times
\sum B_d e^{d\cdot \tau}\cr
(**)~~~\int_X\Omega^{-1}~
\overline{HG[B](t+g(e^t))}~ HG[B](\tau+g(e^\tau))
&=\int_X\Omega^{-1}e^{H_\zeta} 
e^{H\cdot(\bar g(qe^{\zeta\alpha})-g(q))/\alpha}\cr
&\times\sum \bar B_d e^{d\cdot(\tau+g(q))} 
e^{d\cdot(\bar g(qe^{\zeta\alpha})-g(q))+d\cdot\zeta\alpha}\cr
&\times\sum B_d e^{d\cdot(\tau+g(q))}.
}
$$
This shows that the series (**) can be obtained from
(*) by the replacements $\tau\mapsto\tau+g(q)$,
$\zeta\mapsto\zeta+(\bar g(q e^{\zeta\alpha})-g(q))/\alpha$.
Thus combining \dumb~ and \THREE, we get
\eqn\dumbII{
\int_X\Omega^{-1}~
\overline{HG[\tilde B](t)}~ HG[\tilde B](\tau)
=\sum_{d,m} 
q^d e^{d\cdot g(q)}~
\left(\zeta+(\bar g(q e^{\zeta\alpha})-g(q))/\alpha\right)^m
~B'_{d,m}(\alpha).
}
Now write $g=g_+ +g_-$ with $\bar g_\pm=\pm g_\pm$.
Obviously for any $g(q)\in\cR[[q]]$,
$g_+(q e^{\zeta\alpha})-g_+(q)\in \alpha\cdot\cR[[q,\zeta]]$.
Since the involution $\omega\mapsto\bar\omega$ on $\cR$ simply
changes the sign of $\alpha$, the fact that $g_-$ is odd
shows that $g_-(q)\in \alpha\cdot\cR[[q]]$. Likewise for
$g_-(q e^{\zeta\alpha})$. This shows that \dumbII~ 
lies in $\cR[[q,\zeta]]$.
This completes our proof in this case.

Second case: set $\tilde B=\mu_f(B)$, ie.
$$
HG[\tilde B](t)=e^{f/\alpha}~HG[B](t).
$$
Again writing
$f\in \cR[[e^t]]$ as $f=f_+ +f_-$ with
$\bar f_\pm=\pm f_\pm$, we get
$$\eqalign{
\int_X\Omega^{-1}~
\overline{HG[\tilde B](t)}~ HG[\tilde B](\tau)
&=e^{-\bar f(e^t)/\alpha}~
e^{f(e^\tau)/\alpha}~
\int_X\Omega^{-1}~
\overline{HG[B](t)}~ HG[B](\tau)\cr
&=e^{-(f_+(q e^{\zeta\alpha})-f_+(q))/\alpha}~
e^{(f_-(q e^{\zeta\alpha})+f_-(q))/\alpha}\cr
&\times\int_X\Omega^{-1}~
\overline{HG[B](t)}~ HG[B](\tau).
}
$$
The right hand side lies in $\cR[[q,\zeta]]$ as before.
$\Box$

All mirror transformations we will use later
will be of the type $\mu_f,\nu_g$ as above. Moreover,
all Euler data we will encounter will have property
(i) of Theorem \UniquenessTheorem. The transformations
$\mu_f,\nu_g$
clearly preserve this property.

\theorem{Suppose that $A,B$ have property (i) of
Theorem \UniquenessTheorem, and that $A,B$ are
linked.
Suppose that $A$ is an Euler data with $deg_\alpha A_d\leq-2$
for all $d\prec0$,
and that there exists
power series $f\in\cR[[K^\vee]]$, 
$g=(g_1,..,g_m)$, $g_j\in\cR[[K^\vee]]$,
all without constant term, such that
\eqn\dumb{e^{f/\alpha}HG[B](t)=\Omega- 
\Omega {H\cdot (t+g)\over\alpha}
+O(\alpha^{-2})}
when expanded in powers of $\alpha^{-1}$.  Then 
$$HG[A](t+g)=e^{f/\alpha}~HG[B](t).$$
}
\thmlab\Applications
\proof
By Theorem \MirrorTransform, $f,g$ define
two mirror transformations $\mu_f$, $\nu_g$, with
\eqn\dumbII{\eqalign{
HG[\tilde B](t)&=e^{f/\alpha}HG[B](t)\cr
HG[\tilde A](t)&=HG[A](t+g)
}}
where $\tilde B=\mu_f(B)$, $\tilde A=\nu_g(A)$.
Now both $\tilde B,\tilde A$ have property (i)
of Theorem \UniquenessTheorem.
(See remark after Theorem \MirrorTransform.)

Since $deg_\alpha A_d\leq-2$,
$HG[\tilde A](t)$ has the same asymtotic form as
 $HG[\tilde B](t)$  in eqn. \dumb ~ $mod~ O(\alpha^{-2})$.
It follows that
$$e^{H\cdot t/\alpha}~
HG[\tilde A-\tilde B](t)\equiv O(\alpha^{-2}),$$
or equivalently
$deg_\alpha (\tilde A_d-\tilde B_d)
\leq -2$.
Thus $\tilde A,\tilde B$ 
satisfy condition (ii) of Theorem \UniquenessTheorem.
Since $A$ is linked to $B$, it follows that 
$\tilde A$ is linked to $\tilde B$.
By Theorem \UniquenessTheorem,
we conclude that $\tilde A=\tilde B$.
Now our assertion follows from eqns. \dumbII.  $\Box$

\remark{The preceding theorem says that one way to compute
$A$ (or $Q$) is by first finding 
an explicit Euler data $B$ linked to $A$,
and then relate $A$ and $B$ via mirror transformations.}

\newsec{From stable map moduli to Euler data}

Fix an admissible balloon manifold with $c_1(X)\geq0$.
Fix a $T$-equivariant multiplicative class $b_T$.
Its nonequivariant limit is denoted by $b$.
Fix a $T$-equivariant bundle of the form $V=V^+\oplus V^-$,
where $V^\pm$ are respectively the convex/concave bundles.
As before, we write
$$\Omega={b_T(V^+)\over b_T(V^-)}.$$
Let $V_d$ be the bundle induced by $V$ on the 0-pointed
degree $d$
stable map moduli of $X$.
Throughout this section, we denote
$$\eqalign{
Q:~~Q_d&:=\varphi_!(\pi^* b_T(V_d))\cr
K_d&:=\int_{M_{0,0}(d,X)} b(V_d)\cr
\Phi&:=\sum K_d~e^{d\cdot t}\cr
A:~~A_d&:=i_0^*Q_d^v.
}
$$
Note that all these objects depend on the choice of
$b_T$ and $V$, though the notations do not reflect this.

\subsec{The Euler data $Q$}

\theorem{(i) $deg_\alpha A_d\leq-2$.
\item{(ii)} If for each $d$ the class $b_T(V_d)$ has homogeneous
degree the same as the degree of
$M_{0,0}(d,X)$, then in
 the nonequivariant limit we have
$$\eqalign{
\int_X e^{-H\cdot t/\alpha}
A_d&=\alpha^{-3}(2-d\cdot t)K_d\cr
\int_X\left(HG[A](t)-e^{-H\cdot t/\alpha}\Omega\right)&
=\alpha^{-3}(2\Phi-\sum t_i {\partial\Phi\over\partial t_i}).
}$$
}
\thmlab\KdTheorem
\proof
Earlier we have proved that
$$A_d=i_0^*Q_d^v=ev_!\left({\rho^* b_T(V_d)\over
\alpha(\alpha-c_1(L))}\right),$$ where $L=L_d$ is the line 
bundle on ${\cal M}_{0, 1}(d, X)$ whose fiber at a 
point $(f, C;x)$ is the tangent line at $x$.

Assertion (i) now follows immediately from this formula.
The second equality in assertion (ii) follows from
the first equality. 
By the above formula again,
$$\eqalign{
I&:=\int_X e^{-H\cdot t/\alpha}A_d\cr
&=\int_{M_{0,1}(d,X)} e^{-ev^*H\cdot t/\alpha}
{\rho^* b(V_d)\over
\alpha(\alpha-c_1(L))}\cr
&=\int_{M_{0,0}(d,X)}  b(V_d)~
\rho_!\left({e^{-ev^*H\cdot t/\alpha}\over
\alpha(\alpha-c_1(L))}\right).
}$$
Now $b(V_d)$ has homogeneous degree the same as the dimension
$M_{0,0}(d,X)$.
The second factor in the last integrand contributes
a scalar factor given by integration
over a generic fiber $E$ (which is a $\P^1$)
of $\rho$.
So we pick out the degree 1 term in
${e^{-ev^*H\cdot t/\alpha}
\over \alpha(\alpha-c_1(L))}$, which is just
${-ev^*H\cdot t\over \alpha^3}+{c_1(L)\over \alpha^3}$.
Restricting to the generic fiber $E$,
say over $(f,C)\in M_{0,0}(d,X)$, the
evaluation map $ev$ is equal to $f$, which is
a degree $d$ map $E\cong\P^1\ra X$.
It follows that
$$\int_E ev^*H=d.$$
Moreover, since $c_1(L)$ restricted to $E$ is just the
first Chern class of the tangent bundle to $E$, it follows that
$$\int_E c_1(L)=2.$$
So we have
$$I=(-{d\cdot t\over \alpha^3}+{2\over \alpha^3})K_d.~~~~\Box$$

\theorem{
More generally suppose $b_T$ is an equivariant
multiplicative class of the form
$$b_T(V)=x^r+x^{r-1}b_1(V)+\cdots+b_r(V),~~rk~V=r$$
where $x$ is a formal variable, $b_i$ is a characteristic
class of degree $i$.
Suppose $s:=rk~V_d-dim~
M_{0,0}(d,X)\geq0$ is independent of $d\succ0$. Then
$$\eqalign{
{1\over s!}\left({d\over dx}\right)^s|_{x=0}\int_X e^{-H\cdot t/\alpha}
A_d&=\alpha^{-3}x^{-s}(2-d\cdot t)K_d\cr
{1\over s!}\left({d\over dx}\right)^s|_{x=0}
\int_X\left(HG[A](t)-e^{-H\cdot t/\alpha}\Omega\right)&
=\alpha^{-3}x^{-s}(2\Phi-\sum t_i {\partial\Phi\over\partial t_i}).
}$$
}
\thmlab\KdTheoremII
\proof
The proof is entirely analogous to (ii) above. $\Box$

\subsec{Linking theorem for $A$}

Now consider a  mixed bundle $V=V^+\oplus V^-$
on $X$.
Fixed a choice of
equivariant multiplicative class $b_T$.
We assume that $V$
has the following property: there exists
nontrivial $T$-equivariant line bundles
$L_1^+,..,L_{N^+}^+;L^-_1,..,L^-_{N^-}$ on $X$ with
$c_1(L_i^+)\geq0$ and $c_1(L_j^-)<0$, such that
for any balloon $pq\cong\P^1$ in $X$ we have
$$V^\pm|_{pq}=\oplus_{i=1}^{N^\pm}L_i^\pm|_{pq}.$$
Note that $N^\pm=rk~V^\pm$.
We also require that 
$$b_T(V^+)/b_T(V^-)=\prod_i b_T(L^+_i)/\prod_j b_T(L^-_j).$$
In this case we call the list 
$(L_1^+,..,L_{N^+}^+;L^-_1,..,L^-_{N^-})$ the
splitting type of $V$. Note that $V$ is not assumed to
split over $X$. Given such a bundle $V$ and a choice of
multiplicative class $b_T$, we obtain an Euler data 
$Q:~Q_d=\varphi_!(\pi^*b_T(V_d))$ (or $A$) as before.

\theorem{Let $b_T=e_T$ be the equivariant Euler class.
Let $pq$ be a balloon, $d=\delta[pq]\succ0$,
and $\lambda$ be the weight on the tangent line $T_q(pq)$. Then
at $\alpha=
\lambda/\delta$, we have
$$j^*_0(Q_d)|_q
=\prod_i\prod_{k=0}^{\bra c_1(L^+_i),d\ket}
\left(c_1(L^+_i)|_q-k\lambda/\delta\right)\times
\prod_j\prod_{k=1}^{-\bra c_1(L^-_j),d\ket-1}
\left(c_1(L^-_j)|_q+k\lambda/\delta\right).
$$
In particular $Q$
is linked to
$$P:~~P_d=
\prod_i\prod_{k=0}^{\bra c_1(L^+_i),d\ket}
(\hat L^+_i-k\alpha)\times
\prod_j\prod_{k=1}^{-\bra c_1(L^-_j),d\ket-1}
(\hat L^-_j+k\alpha).
$$
}
\proof We first consider one positive line bundle $L$. As in \LLY, we
consider a point $(f, C)\in M_d(X)$ where $f$ is $\delta$-cover from
$C=\P^1$ to the balloon $pq\simeq \P^1$. For
$\alpha=\lambda/\delta$, this map can be written as 

$$f:\ C\rightarrow \P^1\times pq\subset\P^1\times X$$where the
second map is the inclusion.  In terms of coordinates we can write the
first map as  

$$f: \ [w_0, w_1]\rightarrow [w_1,w_0]\times[w_0^\delta, w_1^\delta].$$ 
Note that the $T$-action induces standard rotation on $pq\simeq \P^1$
with the weights $\lambda_1, \lambda_2$ and
$\lambda=\lambda_1-\lambda_2$. It is now easy to see that 
this point $(f, C)$ is fixed by the
subgroup of $G$ with
$\alpha=\lambda/\delta$.
 On the other hand as argued in \LLY,
$(\pi_2\circ f, C)$ is then a smooth fixed point in $M_{0, 0}(d,
X)$ under the $T$-action. The restriction $j_0^*Q_d|_p$ with
$\alpha=\lambda/\delta$ is equal to the value of 
$e_T(\cU_d)$ at $(f, C)$. This, in turn,
is equal to the restriction of $e_T(V_d)$
at $ (\pi_2\circ f, C)$ in $M_{0, 0}(d, X)$.

 Assume the restriction of $L$ to $pq\simeq \P^1$ is $\cO(l)$ with
$l=\bra c_1(L), [pq]\ket$. 
We compute that the equivariant Euler
class restricted to this point  $(\pi_2\circ f, C)$.
As in \LLY,  we get

$$e_T(U_d)=\prod_{m=0}^{l\delta}(l\lambda_1-m\frac{\lambda}{\delta}).$$
Also note that $c_1(L)(p)=l\lambda_1$ and $d=\delta[pq]$, this implies
that $Q_d=\varphi_!(\pi^*e_T(V_d))$ is linked to 

$$P_d=\prod_{m=0}^{\bra c_1(L),d\ket}(c_1(\hat{L})-m\alpha).$$

Similarly for a concave line bundle $L$, if its restriction to the
balloon $pq$ is $\cO(-l)$ with $-l=\bra c_1(L), [pq]\ket$, then 

$$e_T(U_d)=\prod_{m=-}^{l\delta-1}(-l\lambda_1+m\frac{\lambda}{\delta})$$
which implies the formula that in this case $Q_d$ is linked to 

$$P_d=\prod_{m=1}^{-\bra c_1(L),d\ket -1}(c_1(\hat{L})+m\alpha).$$

The general case is just a product of these cases. 
$\Box$

Similarly we can prove the following formula for the Chern polynomial.

\theorem{Let $b_T=c_T$ be the equivariant Chern polynomial.
Let $pq$ be a balloon, $d=\delta[pq]\succ0$,
and $\lambda$ be the weight on the tangent line $T_q(pq)$. Then
at $\alpha=
\lambda/\delta$, we have
$$j^*_0(Q_d)|_q
=\prod_i\prod_{k=0}^{\bra c_1(L^+_i),d\ket}
\left(x+c_1(L^+_i)|_q-k\lambda/\delta\right)\times
\prod_j\prod_{k=1}^{-\bra c_1(L^-_j),d\ket-1}
\left(x+c_1(L^-_j)|_q+k\lambda/\delta\right).
$$
In particular $Q$
is linked to
$$P:~~P_d=
\prod_i\prod_{k=0}^{\bra c_1(L^+_i),d\ket}
(x+\hat L^+_i-k\alpha)\times
\prod_j\prod_{k=1}^{-\bra c_1(L^-_j),d\ket-1}
(x+\hat L^-_j+k\alpha).
$$
}

\newsec{Applications}

\subsec{Toric manifolds} 

We call a toric manifold $X$ {\it reflexive} if
its defining fan satisfies the following combinatorial
condition: the convex hull of
the primitive generators of the 1-cones
in the fan is a reflexive polytope.
It has been shown
\Batyrev\Roan~that a pair of polar reflexive polytopes
gives rise to a pair of mirror (in the sense of
Hodge numbers) Calabi-Yau varieties, by taking
anti-canonical hypersurfaces in the corresponding
reflexive toric manifolds.
It has been conjectured that \BB~ a similar
statement holds for complete intersections in toric manifolds.
It is known that \HLY~
a toric manifold $X$ is reflexive iff
$c_1(X)\geq0$. We shall assume that $X$ is reflexive.
Recall that for a (convex) toric manifold $X$, we have
$$e_G(X_0/W_d)=
\prod_a
\prod_{k=1}^{\bra D_a,d\ket}(D_a-k\alpha)
$$
where each $D_a$ is the $T$-equivariant first Chern classes
of the line bundles corresponding to a $T$-invariant hypersurfaces
in $X$. 

\subsec{Chern polynomials for mixed bundles}

To proceed, we make two further choices:
let $b_T$ be the $T$-equivariant
Chern polynomial $c_T$,
and let
$V=V^+\oplus V^-$ be a mixed bundle with splitting type
$(L^+_1,..,L^+_{N^+}; L^-_1,..,L^-_{N^-})$.
Here the $L$'s are 
$T$-equivariant line bundles on $X$ with
$$\eqalign{&c_1(L^+_i)\geq0,~~c_1(L^-_j)<0,\cr
&\Omega:=c_T(V^+)/c_T(V^-)=\prod_i(x+c_1(L^+_i))
/\prod_j(x+c_1(L^-_j))\cr
&\sum_i c_1(L_i^+)-\sum_j c_1(L^-_j)=c_1(X).
}$$
From this, we get an $\Omega$-Euler data $Q:~Q_d=\varphi_!
(\pi^*c_T(V_d))$ as before.
By the Linking Theorem, $Q$ is linked to the Euler data
$$P:~~P_d=
\prod_i\prod_{k=0}^{\bra c_1(L^+_i),d\ket}
(x+\hat L^+_i-k\alpha)\times
\prod_j\prod_{k=1}^{-\bra c_1(L^-_j),d\ket-1}
(x+\hat L^-_j+k\alpha).
$$
As before, we set
$$B:~~B_d=i_0^*P_d^v,~~~~~A:~~A_d=i_0^*Q_d^v.$$
We consider three separate cases.
We will be using the elementary formula
\eqn\elem{\prod_{k=1}^M({\omega\over\alpha}-k)
\equiv(-1)^M M!(1-{\omega\over\alpha}\sum_{k=1}^M{1\over k})
}
where $``\equiv''$ here means equal $mod~O(\alpha^{-2})$,
to compute the leading terms of
\eqn\Ratio{\eqalign{
B_d&=\prod_i\prod_{k=0}^{\bra c_1(L^+_i),d\ket}(x+c_1(L^+_i)-k\alpha)
\times
\prod_j\prod_{k=1}^{-\bra c_1(L^-_i),d\ket-1}(x+c_1(L^-_i)+k\alpha)
\cr &\times
{1\over
\prod_a\prod_{k=1}^{\bra D_a,d\ket}(D_a-k\alpha)}\cr
&=\Omega~c_T(V^-)~\alpha^{-N^-}~
\prod_i\prod_{k=1}^{\bra c_1(L^+_i),d\ket}
({x+c_1(L^+_i)\over\alpha}-k)
\times
\prod_j\prod_{k=1}^{-\bra c_1(L^-_i),d\ket-1}
({x+c_1(L^-_i)\over\alpha}+k)
\cr &\times
{1\over
\prod_a\prod_{k=1}^{\bra D_a,d\ket}
({D_a\over\alpha}-k)}.
} }

First suppose that $rk~V^-=N^-\geq2$.
In this
case we have 
$$deg_\alpha B_d
=-rk~V^-\leq-2$$
and hence
$$HG[B](t)\equiv\Omega-\Omega{H\cdot t\over\alpha}.$$
By Theorem \Applications, we conclude that $A=B$ and $Q=P$.
This completes the computation of $A$ and $Q$ in this case.

Now consider the case $rk~V^-=N^-=1$, hence $V^-$ is a line bundle.
In this case we have
$$\eqalign{
B_d&\equiv\alpha^{-1}~
\Omega~(x+c_1(V^-))~(-1)^{\bra c_1(V^-),d\ket}
(-\bra c_1(V^-),d\ket-1)!
{\prod_i\bra c_1(L^+_i),d\ket!\over \prod_a\bra D_a,d\ket!}\cr
&=:\alpha^{-1}~\Omega~
(\sum_i H_i\phi_{d,i}+\psi_d)
}$$
where the $\phi_{d,i}\in\Q$, $\psi_d\in\Q[\cT^*,x]$, are determined
uniquely by the writing $c_1(V^-)\in H^2_T(X)$
in the last equality,
according to the decomposition
$H^2_T(X)
=\oplus_{i=1}^m\Q H_i\oplus\cT^*$. Hence we get
$$
e^{-H\cdot t/\alpha}~B_d
\equiv\Omega~
(\alpha^{-1}H\cdot \phi_d+\alpha^{-1}\psi_d).
$$
Summing over $d\in K^\vee$, we get
$$\eqalign{
HG[B](t)&\equiv
\Omega(1-\alpha^{-1} H\cdot(t+F)+\alpha^{-1}G)\cr
F&:=-\sum{}\phi_d ~e^{d\cdot t}\cr
G&:=\sum{}\psi_d ~e^{d\cdot t}.
}
$$
From this we get
$$\eqalign{
 e^{-G/\alpha}~HG[B](t)
&\equiv\Omega-\Omega{H\cdot(t+F)\over\alpha}
}$$
By Theorem \Applications, we conclude that
\eqn\dumb{~e^{-G/\alpha}HG[B](t)
=HG[A](t+F).}
This completes our computation of $A$ and $Q$ in this case.
 
Recall that
$$\eqalign{
dim~ M_{0,0}(d,X)&=\bra c_1(X),d\ket+n-3\cr
rk~V_d&=\sum_i\bra c_1(L^+_i),d\ket-\sum_j\bra c_1(L^-_j),d\ket
+N^+-N_-\cr
&=\bra c_1(X),d\ket+rk~V^+-rk~V^-.}
$$
To applied Theorem \KdTheoremII, we assume that 
$rk~V^+-rk~V^-\geq n-3$,
and we determine all $K_d$ immediately.
Explicitly
(in the
nonequivariant limit $\cT^*\ra0$):
\eqn\MixedMirrorFormula{
{1\over s!}\left({d\over dx}\right)^s|_{x=0}
\int_X\left( e^{-G/\alpha} HG[B](t)-
e^{-H\cdot\tilde t /\alpha}\Omega\right)
=\alpha^{-3}x^{-s}(2\Phi(\tilde t)
-\sum_i\tilde t_i{\partial\Phi(\tilde t)\over\partial\tilde t_i}).
}
where $s:=rk~V^+-rk~V^--(n-3)$, $\tilde t:=t+F(t)$.
Note that this same formula applies also 
when $rk~V^-\geq2$, whereby we put $G,F=0$.

We now consider the case when $V$ is purely convex: $N^-=0$.

\subsec{Convex bundle}

We will denote the $L^+_i$ simply by $L_i$.
Using formulas \elem ~ and \Ratio, we get
$$\eqalign{
B_d&\equiv\Omega
{\prod_i\bra c_1(L_i),d\ket!\over
\prod_a\bra D_a,d\ket!}
(1+\alpha^{-1}\sum_a D_a
\sum_{k=1}^{\bra D_a,d\ket}{1\over k}
-\alpha^{-1}\sum_i(x+ c_1(L_i))
\sum_{k=1}^{\bra c_1(L_i),d\ket}{1\over k}) \cr
&=:\Omega
\lambda_d+\alpha^{-1}\sum_i H_i\phi_{d,i}+\alpha^{-1}\psi_d
}$$
Here
 the $\lambda_d,\phi_{d,i}\in\Q$, 
 $\psi_d\in\Q[\cT^*,x]$
 are determined
uniquely by the writing each $D_a,c_1(L_i)\in H^2_T(X)$
in the last equality,
according to the decomposition
$H^2_T(X)
=\oplus_{i=1}^m\Q H_i\oplus\cT^*$. Since $e^{-H\cdot t/\alpha}
\equiv1-\alpha^{-1}H\cdot t$, we get
$$
e^{-H\cdot t/\alpha}~B_d
\equiv\Omega(
\lambda_d-\alpha^{-1}H\cdot (\lambda_dt-\phi_d)+\alpha^{-1}\psi_d).
$$

Summing over $d\in K^\vee$, we get
$$\eqalign{
HG[B](t)&\equiv
\Omega\left(F_0-\alpha^{-1}
H\cdot(F_0 t+F)+\alpha^{-1}G\right)\cr
F_0&:=1+\sum \lambda_d~ e^{d\cdot t}\cr
F&:=-\sum\phi_d ~e^{d\cdot t}\cr
G&:=\sum\psi_d ~e^{d\cdot t}.
}
$$

Put $f:=\alpha~log~F_0-{G\over F_0}$.
Then we get
$$\eqalign{
 e^{f/\alpha}~HG[B](t)
&\equiv\Omega-\Omega{H\cdot(t+{F\over F_0})\over\alpha}
}$$
By Theorem \Applications, we conclude that
\eqn\dumb{e^{f/\alpha}~HG[B](t)
=HG[A](t+{F\over F_0}).}
This completes our computation of $A$ and $Q$ in this case.

Again to apply Theorem \KdTheoremII, we assume that
$rk~V\geq n-3$, and
determine all $K_d$ immediately. 
Explicitly:
\eqn\MirrorFormula{
{1\over s!}\left({d\over dx}\right)^s|_{x=0}
\int_X\left( e^{f/\alpha}~HG[B](t)-
e^{-H\cdot\tilde t /\alpha}\Omega\right)
=\alpha^{-3}x^{-s}(2\Phi(\tilde t)
-\sum_i\tilde t_i{\partial\Phi(\tilde t)\over\partial\tilde t_i}).
}
where $s:=rk~V-(n-3)$, $\tilde t:=t+{F(t)\over F_0(t)}$.

Let us now specialize to the
case $rk~V=n-3$ (ie. $s=0$), and $V=\oplus_i L_i$.
We can then set $x=0$,
so that $b_T=c_T$ becomes the equivariant Euler class $e_T$,
and the $K_d$ is just the intersection numbers for $e(V_d)$.
Then the formula \MirrorFormula~ yields
the general formula derived in \HLYII~
and in \HKTY,
on the basis of the conjectural mirror correspondence.
Note that
$$F_0=\sum
{\prod_i\bra c_1(L_i),d\ket!\over
\prod_a\bra D_a,d\ket!} e^{d\cdot t}$$
is an example of
a hypergeometric function \GKZ.
It has been proved in \HLY~ that $F_0$ is the unique
holomorphic period of Calabi-Yau hypersurfaces near
the so-called large radius limit.
For the purpose of comparison, we should mention that the 
definition of $\Phi$ here differs from the prepotential
in \HLYII\HKTY~ by a degree three polynomial in $\tilde t$, and
the definition of the hypergeometric series $HG[B](t)$
here differs from that denoted by $w_0(x,\rho)$ in \HLYII\HKTY~
by an irrelevant
overall constant factor.

Precursors to the above general formula have been many examples
\HKTY\BatyrevStraten\CDFKM\BKK.
We now specialize to a few numerical examples
which have been frequently studied by both physicists
and mathematicians alike.

\subsec{A complete intersection in $\P^1\times\P^2\times\P^2$}

The complete intersection of degrees $(1,3,0),~
(1,0,3)$ in this 5-dimensional toric manifold $X$ has been studied
in \HKTY~ using mirror symmetry, and in \HSS~
computing some of the intersection  numbers $K_d$ for
the Euler class $b=e$ in terms of
modular forms. 

From our point of view, that complete intersection correspond
to the following choice of convex bundle:
$$V=\cO_1(1)\otimes\cO_2(3)\oplus\cO_1(1)\otimes\cO_3(3)$$
where $\cO_i(l)$ denotes the pullback of $\cO(l)$
from the $i$th factor. The K\"ahler cone of $X$
is abviously generated by the hyperplanes $H_1,H_2,H_3$
from the three factors of $X$, and hence
$K^\vee$ can be identified with the set of
$d=(d_1,d_2,d_3)\in\Z_{\geq0}^3$. We consider
intersection numbers $K_d$ for
the Euler class $b=e$ as before.
Thus we set $\Omega=e(V)=(H_1+3H_2)(H_1+3H_3)$.
The Euler data $P$ we need in eqn. \MirrorFormula~
is given by 
$$\eqalign{
P_d&=\prod_{k=0}^{d_1+3d_2}(\kappa_1+3\kappa_2-k\alpha)
\times\prod_{k=0}^{d_1+3d_3}(\kappa_1+3\kappa_3-k\alpha)\cr
j^*_0(P_d)&=\prod_{k=0}^{d_1+3d_2}(H_1+3H_2-k\alpha)
\times\prod_{k=0}^{d_1+3d_3}(H_1+3H_3-k\alpha).}
$$

The linear sigma model is $W_d=N_d(\P(n))= N_{d_1, 1}\times N_{d_2, 2}\times N_{d_3, 2}$. The equivariant Euler class, after taking nonequivariant limit with
 respect to the $T$ action, is given by  

$$e_G(X_0/W_d)=\prod_{m=1}^{d_1}(H_1-m\alpha)^2\prod_{m=1}^{d_2}(H_2-m\alpha)^3\prod_{m=1}^{d_3}(H_3-m\alpha)^3.$$Now we can easily write down the hypregeometric series and all the $K_d$ can be computed by our formula
\MirrorFormula~ at once using the obvious intersection form on $X$,
given by the relations:
$$\int_X H_1 H_2^2 H_3^2=1,~~H_1^2=1,~H_2^3=1,~ H_3^3=1.$$

Once we have the hypergeometric series, the corresponding Picard-Fuchs equation can be easily written down as given in \HKTY.

\subsec{$V=\cO_1(-2)\otimes\cO_2(-2)$ on $\P^1\times\P^1$}

Here we denote by $\cO_i(l)$ the pullback of
$\cO(l)$ from the $i$th factor of $X=\P^1\times\P^1$.
Our bundle $V$ has $rk~V^+-rk~V^-=n-3=-1$. Thus
we can apply our formula \MixedMirrorFormula~
with $x=0$.
We put $\Omega={1\over H_1H_2}$.
The Euler data $P$ in eqn. \MixedMirrorFormula~
that compute the $K_d$ is now given by:
$$
P_d=\prod_{k=1}^{2d_1-1}(-2\kappa_1+k\alpha)\times
\prod_{k=1}^{2d_2-1}(-2\kappa_2+k\alpha).
$$

The corresponding equivariant Euler class, after taking the nonequivariant limit with respext to the $T$-action is  

$$e_G(X_0/W_d)= \prod_{m=1}^{d_1}(H_1-m\alpha)^2\prod_{m=1}^{d_2}(H_2-m\alpha)^2.$$

Again one can immediately write down the hypergeometric series as well as the 
corresponding Picard-Fuchs equation by using our mirror principle.

\newsec{Generalizations and Concluding Remarks} 

\subsec{A weighted projective space}

Consider the following example: the concave bundle $V=O(-6)$ over $\P_{3,2,1}$, $\Omega=\frac{1}{6H}$. This example will be studied in our subsequent paper by using resolution of singularities. 
This is an example of ``local mirror symmetry''
studied in physics \LMW. The mirror formula there
can derived as a special case of our general result. 
In fact, the Euler data which computes
the $K_d$ in this case is determined by

$$j_0^*P_d= \prod_{m=1}^{6d-1}(-6H+m\alpha).$$
The corresponding equivariant Euler class, after taking nonequivariant limit with respect to the $T$ action, is: 

$$e_G(X_0/W_d)= \prod_{m=1}^{d}(H-m\alpha)\prod_{m=1}^{2d}(2H-m\alpha)\prod_{m=1}^{3d}(3H-m\alpha).$$

The corresponding hypergeometric series and Picard-Fuchs 
equation can be immediately written down.
It turns out that the hypergeometric series gives the periods of 
a meromorphic 1-form for a family of elliptic curves \LMW.

\subsec{General projective balloon manifolds}

Let $X$ be a projective manifold embedded in
$\P(n)$,
 with a system of homogeneous polynomial defining
equations 
$P(z^1, \cdots, z^n)=0$, where $z^j=(z^j_1, \cdots, z^j_{n_j})$. 
For each $P$, 
by taking the coefficients of each monomial 
$w_0^aw_1^b$ in $P(f^1, \cdots,
f^n)=0$, where $f^j=[f^j_1(w_0, w_1),\cdots,
f^j_{n_j}(w_0, w_1)]$ for $j=1, \cdots, k$ is the tuple of polynomials that define the coordinates of $N_d(\P(n))$, we get 
several equations of the same degree as $P$. These equations together define a projective variety, which 
we denote by $N_d(X)$, in $N_d(\P(n))$.

As discussed earlier, we see that the
$S^1$ fixed point components in $N_d(X)$ are given by the $X_r$'s which are copies of $X$. We do not know whether the localization formula holds on $N_d(X)$. The localization formula holds 
if the fixed point components embedded into 
$W_d$ as local complete intersection subvarieties. It is likely that this is the case for any convex projective manifold. If this
 is true, then we can take $N_d(X)$ to be the linear sigma model $W_d$ for $X$. Then our mirror principle may apply readily 
to compute multiplicative characteristic numbers on 
$M_{0,0}(d,X)$ in terms of the hypergeometric series. 

\subsec{A General Mirror Formula}

Many of our results so far are proved for projective manifolds
without $T$-action. Here we first discuss a formula for computing the numbers
$$K_d=\int_{M_{0,0}(d,X)}b(V_d)$$
for a general convex projective $n$-fold $X$ {\it without $T$-action}.
For simplicity, let's focus on the case when the multiplicative
class $b$ is the Chern polynomial $c$, and $V$ is a direct
sum of line bundles on $X$. There is a similar
formulation in the general case.  We fix a projective
embedding $X\ra\P(n)$, as before.
Note that the map $\varphi:M_d(X)\ra N_d(\P(n))$ is now
only $S^1$-equivariant. 
Recall that
the subvariety $W_d:=\varphi(M_d(X))\subset N_d(\P(n))$
contains as
$S^1$ fixed point components copies of $X$: $X_r$, $0\preceq r\preceq
d$. We assume that the localization formula holds on it. 

We denote by $e_{S^1}(X_0/W_d)$ the equivariant Euler class of the normal
bundle of $X_0$ in $W_d$.
Let
$$V=V^+\oplus V^-,~~~V^+:=\oplus L_i^+,~~~V^-:= \oplus L_j^-$$
satisfying
$c_1(V^+)-c_1(V^-)=c_1(X)$ and
$rk(V^+)-rk(V^-)-(n-3)\geq0$,
where the $L_i^\pm$ are respectively convex/concave line
bundles on $X$. Let
$$\eqalign{
\Omega&=B_0:=c(V^+)/c(V^-)=\prod_i(x+c_1(L^+_i))
/\prod_j(x+c_1(L^-_j))\cr
B_d&:={1\over e_{S^1}(X_0/W_d)}\times
\prod_i\prod_{k=0}^{\bra c_1(L^+_i),d\ket}
(x+c_1( L^+_i)-k\alpha)\times
\prod_j\prod_{k=1}^{-\bra c_1(L^-_j),d\ket-1}
(x+c_1( L^-_j)+k\alpha).\cr
HG[B](t)&:=\sum B_d e^{d\cdot t}\cr
\Phi(t)&:=\sum K_d e^{d\cdot t}.
}
$$

\conjecture{There exist unique power series $G(t),F(t)$ such that
the following formula holds:
$$
{1\over s!}\left({d\over dx}\right)^s|_{x=0}
\int_X\left( e^{-G/\alpha} HG[B](t)-
e^{-H\cdot\tilde t /\alpha}\Omega\right)
=\alpha^{-3}x^{-s}(2\Phi(\tilde t)
-\sum_i\tilde t_i{\partial\Phi(\tilde t)\over\partial\tilde t_i}).
$$
where $s:=rk~V^+-rk~V^--(n-3)$, $\tilde t:=t+F(t)$.
Moreover $G,F$ are determined by the condition that
the integrand on the left hand side is of order $O(\alpha^{-2})$.
}

\subsec{Formulas without $T$-action }

One of our key ingredient, the functorial localization formula 
plays an important role in relating
 the data on $M_d(X)$ and those on $W_d$. 
It turns out that similar formula holds in $K$-theory. 
It holds even when $X$ has no group action. This
indicates that our method may be extended to compute 
$K$-theory multiplicative type characteristic classes 
on $M_d(X)$ (and ultimately on $M_{0,0}(d,X)$),
in terms of certain $q$-hypergeometric series, even
for projective manifolds without group action.   

We now write down the  relevant localization formulas for 
convex $X$ without torus action, both in cohomology and in $K$-theory. 
The notations and proofs are basically the same as before. 
Given a manifold $X$, let's assume that there is
a linear sigma model $W_d$. 

\lemma{ For any equivariant cohomology class $\omega$ on $M_d(X)$, the following equality holds on $X_r$ for any $0\preceq r\preceq d$:
$$\frac{j_r^*\varphi_!(\omega)}{e_{S^1}(X_r/W_d)}=e_![\frac{i_r^*\omega}{e_{S^1}(F_r/M_d(X))}].$$
}
Here $e_{S^1}(\cdot)$ denotes the $S^1$-equivariant 
Euler class. As in the cases we have studied earlier,
the left hand side of the above formula indicates
that when $V=L$ is a line bundle,
 we should compare the Euler data 
$Q_d=\varphi_!\pi^*e(V_d)$, to the Euler data given by
$$P_d=\prod_{m=0}^{\bra c_1(L), d\ket}(c_1(L)-m\alpha).$$

What is left is to develope  uniqueness and mirror transformations, 
which we are unable to achieve at this moment, though 
they can be easily axiomized.

 Now let us look at $K$-theory formula, which can be proved by using equivariant localization in K-theory. First following the same idea, we get the explicit formula as follows: given any equivariant element $V$ in $K(W_d)$, we have  

$$V=\sum_r{j_r}_!\frac{j_r^*V}{E_G(X_r/W_d)}$$ 
where $E_G(X_r/W_d)$ is the equivariant Euler class of the normal bundle of $X_r$ in $N_d(X)$. 
Here the push-forward and pull-backs by $j_r$ denote the corresponding operations in $K$-theory. By taking $V=1$, we get

 $$e_![\frac{1}{E_G(F_r/M_d(X))}]=\frac{j_r^* \varphi_!(1)}{E_G(X_r/W_d)}.$$

 Second, we have the following lemma; 
 If $\varphi_!(1) =1$, which is the case if $X=\P^n$, this formula
 gives explicit formulas for some $K$-theoretic characteristic numbers of
 the moduli spaces.

\lemma{Given any equivariant element $V$ in $K_G(M_d(X))$, then we have formula $$e_![\frac{i_r^*V}{E_G(F_r/M_d(X))}]=\frac{j_r^*(\varphi_!V)}{E_G(X_r/W_d)}$$ where $E_G(\cdot)$ denotes the equivariant Euler class of the corresponding normal bundle in K-group.}

 In particular we have explicit expressions from the decomposition of the normal bundles:

$$E_G(F_r/M_d(X))=(1-e^{\alpha})(1-e^{-\alpha})(1-e^\alpha L_r)(1-e^{-\alpha} L_{d-r})$$
and similarly 

$$E_G(F_0/M_d(X))=(1-e^{\alpha})(1-e^\alpha L_d),\ E_G(F_d/M_d(X))
=(1-e^{-\alpha})(1-e^{-\alpha} L_d).$$

For a toric manifold $X$, we also have the explicit class in $K$-group,

$$E_G(X_0/W_d)=\prod_a
\prod_{m=1}^{\bra D_a,d\ket }(1-e^{m\alpha}[D_a])$$
where $[D_a]$ are the equivariant line bundle
corresponding to the $T$ divisors $D_a$.

If $V$ is a multipicative type K-theory characteristic class, 
then we can develope a similar theory of Euler data and
uniqueness. 
These result can also be extended to the
nonconvex case without a hitch.

\footatend\vfill\supereject\immediate\closeout\rfile\writestoppt
\baselineskip=14pt\centerline{{\bf References}}\bigskip{\frenchspacing%
\parindent=20pt\escapechar=` \input refs.tmp\vfill\eject}\nonfrenchspacing

\end